\documentclass[sn-mathphys-num]{sn-jnl}

\usepackage{graphicx}%
\usepackage{multirow}%
\usepackage{amsmath,amssymb,amsfonts}%
\usepackage{amsthm}%
\usepackage{mathrsfs}%
\usepackage[title]{appendix}%
\usepackage{xcolor}%
\usepackage{textcomp}%
\usepackage{manyfoot}%
\usepackage{booktabs}%
\usepackage{algorithm}%
\usepackage{algorithmicx}%
\usepackage{algpseudocode}%
\usepackage{listings}%
\usepackage{enumitem}

\usepackage{amsmath} 
\usepackage{amssymb} 
\usepackage{amsthm} 
\usepackage{hyperref} 

\usepackage{graphicx} 
\usepackage{float}    

\theoremstyle{thmstyleone}%
\newtheorem{theorem}{Theorem}
\newtheorem{proposition}[theorem]{Proposition}%

\theoremstyle{thmstyletwo}%
\newtheorem{example}{Example}%
\newtheorem{remark}{Remark}%

\theoremstyle{thmstylethree}%
\newtheorem{definition}{Definition}%
\newtheorem{claim}{Claim}%
\newtheorem{lemma}{Lemma}%
\newtheorem{corollary}[theorem]{Corollary}

\begin{document}

\title[Pseudo-timelike loops in signature changing manifolds]{Pseudo-timelike loops in signature changing semi-Riemannian manifolds with a transverse radical}

\author*[1,2]{\fnm{N. E.} \sur{Rieger}}\email{n.rieger@math.uzh.ch, n.rieger@yale.edu}

\author[3,4]{\fnm{W.} \sur{Hasse}}

\affil*[1]{\orgdiv{Department of Mathematics}, \orgname{University of Zurich}, \orgaddress{\street{Winterthurerstrasse
190}, \city{Zurich}, \postcode{8057}, \country{Switzerland}}}

\affil[2]{{Current address:} \orgdiv{Mathematics Department}, \orgname{Yale University}, \orgaddress{\street{219 Prospect Street}, \city{New Haven}, \postcode{06520}, \state{CT}, \country{USA}}}

\affil[3]{\orgdiv{Institute for Physics and Astronomy}, \orgname{Technical University
Berlin}, \orgaddress{\street{Hardenbergstr. 36}, \city{Berlin}, \postcode{10623}, \country{Germany}}}

\affil[4]{\orgname{Wilhelm Foerster
Observatory Berlin}, \orgaddress{\street{Munsterdamm 90}, \city{Berlin}, \postcode{12169}, \country{Germany}}}


\abstract{In 1983, Hartle and Hawking proposed the no-boundary proposal, 
suggesting that the universe has no beginning in the sense of a spacetime singularity or boundary.
Nevertheless, there is an origin of time. Mathematically, this involves signature-type changing manifolds 
in which a Riemannian region smoothly transitions to a Lorentzian region across the hypersurface $\mathcal{H}$ where time begins. 

We develop a coherent framework for signature changing manifolds 
with a degenerate yet smooth metric. Established Lorentzian tools and results 
are then adapted to this setting, and new definitions are introduced that carry unforeseen causal
implications. A noteworthy consequence is the presence of locally
time-reversing loops through every point on the hypersurface.
Imposing global hyperbolicity on the Lorentzian region, we prove that for every point $p \in M$ 
there exists a pseudo-timelike loop self-intersecting at $p$. Equivalently, $M$ always admits a closed 
pseudo-timelike path around which the time direction reverses, preventing 
any consistent distinction between future- and past-directed vectors. To an 
observer near $\mathcal{H}$, such loops may appear as the creation of a 
particle-antiparticle pair at two distinct points.}




\maketitle

\section{Introduction}\label{sec1}

According to popular ideas about quantum cosmology, classical cosmological
models contain an initial Riemannian region of Euclidean signature
joined to a final semi-Riemannian region with the usual Lorentzian
signature~\cite{Ellis - Change of signature in classical relativity,Ellis - Covariant change of signature in classical relativity}.
In 1983 Hartle and Hawking~\cite{Hartle Hawking - Wave function of the Universe}
introduced a conceptually intriguing idea involving signature-type
change, which led to the no-boundary proposal for the initial conditions
of the universe.\footnote{Although the Hartle-Hawking proposal does not explicitly address signature change, it naturally motivates the concept: the universe ``tunnels'' into existence from a state with no time, a Euclidean-signature spacetime without a past boundary, to a Lorentzian spacetime. This process effectively eliminates the classical Big Bang singularity. The resulting transition and its classical implications have been discussed in~\cite{Ellis - Change of signature in classical relativity,Vilenkin}.} In this proposal, the universe has no initial singularity or boundary.\footnote{Although singularities can be considered points where curves terminate
at finite parameter values, providing a general definition remains
difficult~\cite{Geroch - What is a singularity in general relativity}.}  While such a spacetime is singularity-free and thus has no distinct beginning, it still possesses an origin of time~\cite{Ellis - Change of signature in classical relativity}.

\textbf{~}\\Since a signature-type changing metric is necessarily either degenerate
or discontinuous at the locus of signature change~\cite{Dray - Gravity and Signature Change},
we will allow for the metric to become degenerate. Hence, in the present
article we will discuss singular semi-Riemannian manifolds for which
the metric constitutes a smooth $(0,2)$-tensor field that is degenerate
at a subset $\mathcal{H}\subset M$, where the bilinear type of the
metric changes upon crossing $\mathcal{H}$.

~

Although the compatibility of the Riemannian and Lorentzian domains
is assumed to be established, insofar as the metric should be smooth
on the interface $\mathcal{H}$, the behavior of curves as they cross
this interface still requires further study. Moreover, in a manifold
where the signature changes from $(+,+,\ldots,+)$ to $(-,+,\ldots,+)$,
the conventional concept of timelike (or spacelike) curves does not
exist anymore. This gives rise to a new notion of curves called \textit{pseudo-timelike}
and \textit{pseudo-spacelike curves}. In order to define these curves
we make a detour to draw upon the concept of the generalized affine
parameter which we use as a tool to distinguish genuine pseudo-timelike
(and pseudo-spacelike, respectively) curves from curves that asymptotically
become lightlike as they approach the hypersurface of signature change. 

~

We endeavor to adapt well-established Lorentzian tools and results
to the signature changing setting, as far as possible. This task proves
to be less straightforward than anticipated, necessitating the introduction
of new definitions with unexpected causal implications, reaching a
critical juncture in our exploration. We draw upon the definition
of pseudo-time orientability and the given absolute time function
to decide whether a pseudo-timelike curve is future-directed. This
establishes the definition for the pseudo-chronological past (and
pseudo-chronological future) of an event. 

~

In this article we show that for signature-type change of the delineated
type, all these considerations lead to a surprising theorem revealing
the non-well-behaved nature of these manifolds. In a sufficiently
small region near the junction of signature change $\mathcal{H}$,
transverse signature-type changing manifolds with a transverse radical
exhibit local anomalies: Specifically, each point on the junction
facilitates a closed time-reversing loop, challenging conventional
notions of temporal consistency.\footnote{\textcolor{black}{In more informal terms, in general relativity, a
closed timelike curve is a smooth, timelike loop where, at every intersection
point, the direction of movement is consistently the same. In contrast,
a loop is a broader concept where a timelike curve loops back on itself,
but the direction of movement at the intersection points might not
always be the same. This is a more intuitive explanation; for a precise
mathematical definition and its extension to a setting with signature-type
change, see Definition~\ref{Definition Chronology violating curve}.}} 
Or put differently, there always exists a closed pseudo-timelike path in $M$ along which the direction of time reverses, making it impossible to consistently assign future- and past-directed vectors. By imposing global hyperbolicity on the Lorentzian region, the global analog can be proven by showing
that for every point $p\in M$, there exists a pseudo-timelike loop intersecting at $p$. In other words, a transverse, signature-type changing manifold with a transverse radical always admits pseudo-timelike loops. From the viewpoint of an observer in the Lorentzian region near the hypersurface $\mathcal{H}$, such locally pseudo-timelike loops may be interpreted as the creation of a particle-antiparticle pair at two distinct points $\hat{q},q \in \mathcal{H}$.

\subsection{Transverse type-changing singular semi-Riemannian manifold}

Unless otherwise specified, the considered manifolds, denoted as $M$
with dimension $\dim(M)=n$, are assumed to be locally homeomorphic
to $\mathbb{R}^{n}$. Moreover, these manifolds are expected to be
connected, second countable, and Hausdorff. This definition also indicates
that all manifolds have no boundary. Additionally, we will generally
assume that the manifolds under consideration are smooth. Unless stated
otherwise, all related structures and geometric objects (such as curves,
maps, fields, differential forms, etc.) are assumed to be smooth as
well.

~

\begin{definition}
A \textit{singular semi-Riemannian manifold} is a generalization of
a semi-Riemannian manifold. It is a differentiable manifold having
on its tangent bundle a symmetric bilinear form which is allowed to
become degenerate.
\end{definition}

~

\begin{definition}
\label{Definition. signature change} Let $(M,g)$ be a singular semi-Riemannian
manifold and let be $p\in M$. We say that the metric changes its
signature at a point $p\in M$ if any neighborhood of $p$ contains
at least one point $q$ where the metric's signature differs from
that at $p$.
\end{definition}

~

We align with~\cite{Kossowksi + Kriele - Signature type change and absolute time in general relativity}
in requiring that $(M,g)$ be a semi-Riemannian manifold with $\dim M\geq2$,
where $g$ is a smooth, symmetric, degenerate $(0,2)$-tensor on $M$,
and $\mathcal{H}:=\{q\in M\!\!:g\!\!\mid_{q} is\;degenerate\}$.{\small{}
}This means $\mathcal{H}$ is the locus where the rank of $g$ fails
to be maximal. In addition\textcolor{black}, we assume that one connected
component of $M\setminus\mathcal{H}$ is Riemannian, denoted by $M_{R}$,
while all other connected components $(M_{L_{\alpha}})_{\alpha\in I}\subseteq M_{L}\subset M$
are Lorentzian, where $M_{L}:=\underset{\alpha\in I}{\bigcup}M_{L_{\alpha}}$represents
the Lorentzian domain. Furthermore, we assume throughout that the
point set $\mathcal{H}$, where $g$ becomes degenerate is not empty.

~ 

\begin{remark}
This restriction permits transitions of the form \( M_L - M_R - M_L \); however, it forbids \( M_R - M_L - M_R \), which could correspond to transitions involving both an initial and a final singularity. This limitation is of a physical nature, ensuring that the Hartle-Hawking ``no boundary" proposal is satisfied, thereby allowing the Lorentzian part \( M_L \) to be interpreted as our universe. Then within \( M_L \), Einstein's equations hold—these form a hyperbolic system of partial differential equations whose solution (given known matter conditions, such as a vacuum) is determined by a Riemannian hypersurface and the time derivatives of the induced metric on that hypersurface. 

~

On the other hand, hypersurfaces that permit a signature change are highly constrained, as the evolution of the solution for \( g \) would have to degenerate across an entire hypersurface. Consequently, we assume that the Lorentzian part cannot revert to a region of Riemannian nature. For the Riemannian part, this restriction does not apply. First, there is no reason to assume Einstein’s equations for conventional matter in this region. Second, the system would no longer be hyperbolic in nature. From a purely mathematical perspective, this argument is, of course, irrelevant—at the level of the metric \( g \), multiple transitions can certainly be allowed.  
\end{remark}

~\\ Moreover, we impose the following two conditions~\cite{Kossowksi + Kriele - Signature type change and absolute time in general relativity}:

~
\begin{enumerate}
\item We call the metric $g$ a codimension-$1$\textbf{ transverse type-changing
metric} if $d(\det([g_{\mu\nu}]))_{q}\neq0$ for any $q\in\mathcal{H}$
and any local coordinate system $\xi=(x^{0},\ldots,x^{n-1})$ around
$q$. Then we call $(M,g)$ a \textbf{transverse type-changing singular
semi-Riemannian manifold}~\cite{Aguirre+Lafuente, Aguirre - On the Conformal Geometry of Transverse Riemann-Lorentz Manifolds, Kossowksi + Kriele - Signature type change and absolute time in general relativity}.~\\ This
implies that the subset $\mathcal{H}\subset M$ is a smoothly embedded
hypersurface in $M$, and the bilinear type of $g$ changes upon crossing
$\mathcal{H}$. Moreover, at every point $q\in\mathcal{H}$ there
exists a one-dimensional subspace, denoted as the radical $Rad_{q}\subset T_{q}M$,
within the tangent space $T_{q}M$ that is orthogonal to all of $T_{q}M$
at that point.~\\{}
\item \textbf{The radical $Rad_{q}$ is transverse to $\mathcal{H}$ for any $q\in\mathcal{H}$.} 
In cosmological applications, particularly in the context of “no boundary” models, a spacelike surface of signature change is the natural and preferred choice. Thus, our focus is specifically directed toward the concept of a transverse radical, which ensures such a spacelike hypersurface of signature change.
Henceforward, we assume throughout that $(M,g)$ is a singular transverse
type-changing semi-Riemannian manifold with a \textit{transverse radical},
unless explicitly stated otherwise.\footnote{ In the case of a tangential radical, no statements corresponding to the main results of our work can be made, since in $M_{L}$ causal curves approaching the hypersurface tend to become tangent to it rather than intersecting it transversely. This corresponds precisely to the case of a lightlike surface, as in some models of signature change in black hole spacetimes.}
\end{enumerate}

~

\begin{remark}
Recall that the radical at $q\in\mathcal{H}$ is defined as the subspace
$Rad_{q}:=\{w\in T_{q}M:g(w,\centerdot)=0\}$. This means $g(v_{q},\centerdot)=0$
for all $v_{q}\in Rad_{q}$. Note that the radical can be either transverse
or tangent to the hypersurface $\mathcal{H}$. The radical $Rad_{q}$
is called \textit{transverse}~\cite{Kossowski - The Volume BlowUp and Characteristic Classes for Transverse}
if $Rad_{q}$ and $T_{q}\mathcal{H}$ span $T_{q}M$ for any $q\in\mathcal{H}$,
i.e. $Rad_{q}+T_{q}\mathcal{H}=T_{q}M$. This means that $Rad_{q}$
is not a subset of $T_{q}\mathcal{H}$, and obviously, $Rad_{q}$
is not tangent to $\mathcal{H}$ for any $q$.
\end{remark}

~

The following theorem is a direct consequence of the above two conditions and the results in~\cite{Kossowksi + Kriele - Signature type change and absolute time in general relativity,Kriele + Martin - Black Holes ...},
but it is restated here in a more concise and self-contained form for clarity. Moreover, we have refined its statement for improved readability and conceptual transparency:

~
\begin{theorem}
\label{thm:Radical-adapted Gauss-like coordinates} Let $M$ be a
singular semi-Riemannian manifold endowed with a $(0,2)$-tensor field
$g$ and the surface of signature change defined as $\mathcal{H}:=\{q\in M\!\!:g\!\!\mid_{q} is\;degenerate\}$.
Then $(M,g)$ is a transverse, signature-type changing manifold with
a transverse radical if and only if for every $q\in\mathcal{H}$ there
exist a neighborhood $U(q)$ and smooth coordinates $(t,x^{1},\ldots,x^{n-1})$
such that $g=-t(dt)^{2}+g_{ij}(t,x^{1},\ldots,x^{n-1})dx^{i}dx^{j}$,
for $i,j\in\{1,\ldots,n-1\}$.
\end{theorem}

~

In the style of time-orthonormal coordinates in Lorentzian geometry
we denote the coordinates in Theorem~\ref{thm:Radical-adapted Gauss-like coordinates}
as \textbf{radical-adapted Gauss-like coordinates}. It is now possible
to simplify matters by using these coordinates whenever dealing with
a transverse, signature-type changing manifold with a transverse radical.
Notably, signature-type change and the radical-adapted Gauss-like
coordinates imply the existence of an uniquely determined, coordinate
independent, natural\textit{ absolute time function $\mathfrak{h}(t,\mathbf{\hat{x}}):=t$}, with $\mathbf{\hat{x}}=(x^{1},\ldots,x^{n-1})$,
in the neighborhood of the hypersurface~\cite{Kossowksi + Kriele - Signature type change and absolute time in general relativity}.
Then the absolute time function establishes a foliation~\cite{Gibbons + Hartle,Halliwell + Hartle}
in a neighborhood of $\mathcal{H}$, such that $\mathcal{H}$ is a
level surface of that decomposition.

\subsection{Statement of results}

Before presenting the main results we require some new definitions.

~

\begin{definition}[Pseudo-timelike curve]
Let $\gamma\colon [a,b] \to M$ be a continuous and differentiable
curve, with $[a,b]\subset\mathbb{R}$ and $-\infty<a<b<\infty$.  
We call $\gamma=\gamma^{\mu}(u)=x^{\mu}(u)$ in $M$ a 
\emph{pseudo-timelike} (respectively, \emph{pseudo-spacelike}) curve if
\begin{enumerate}
\item $\operatorname{Im}(\gamma) \cap M_{L} \neq \varnothing$, i.e.\ $\gamma$ has image points in the Lorentzian region; and
\item for every generalized affine parametrization (see Definition~\ref{Definition-GAP-(Beem):Definition Generalized Affine Parameter}) of $\gamma$ in $M_{L}$ there exists $\varepsilon>0$ such that  $g(\gamma',\gamma') < -\varepsilon$ (respectively, $g(\gamma',\gamma') > \varepsilon$).

\end{enumerate}
\end{definition}

\textbf{~}\\In simpler terms, we call a curve pseudo-timelike if it is timelike
in the Lorentzian domain $M_{L}$ and does not become asymptotically
lightlike as it approaches the hypersurface where the signature changes.
Consequently, a pseudo-timelike loop is a generalization of a pseudo-timelike
curve that loops back on itself. However, unlike a regular closed
curve where the direction of movement would be the same at every intersection
point, in a pseudo-timelike loop, the direction of movement at the
intersection points is not necessarily the same (see Definition~\ref{Definition Chronology violating curve}).

~ 

{
\renewcommand\thedefinition{}
\begin{definition}[Pseudo-timelike]
\textbf{\label{def:Pseudo-timelike-1}}A vector
field $V$ on a signature-type changing manifold $(M,g)$ is \textit{pseudo-timelike}
if and only if its integral curves are pseudo-timelike (hence, in particluar, $\,V$ is timelike in $M_{L}$).
\end{definition}
}

~

{
\renewcommand\thedefinition{}
\begin{definition}[Pseudo-time orientable]
\textbf{\label{def:Pseudo-time-orientable-1}}
A signature-type changing manifold $(M,g)$ is \textit{pseudo-time
orientable} if and only if the Lorentzian region $M_{L}$ is time
orientable.
\end{definition}
}

~

In a sufficiently small region near the junction of signature change,
transverse, signature-type changing manifolds with a transverse radical
exhibit local anomalies. Specifically, each point on the junction
gives rise to the existence of closed time-reversing loops, challenging
conventional notions of temporal consistency.

~

\begin{theorem}[Local loops]
\label{Theo: Local Loops-1}Let $(M,\tilde{g})$ be a transverse,
signature-type changing, $n$-dimensional ($n\geq2$) manifold with
a transverse radical. Then in each neighborhood of each point ${\color{black}\ensuremath{q\in\mathcal{H}}}$
there always exists a pseudo-timelike loop.
\end{theorem}

~  

The existence of such pseudo-timelike curves locally near the hypersurface
that loop back to themselves, gives naturally rise to the question
whether this type of curves also occur globally. In the global version
a key notion is global hyperbolicity which plays a role in the spirit
of completeness for Riemannian manifolds. By imposing the constraint
of global hyperbolicity on the Lorentzian region, we demonstrate

~

\begin{theorem}[Global loops]
\label{Theo: Global Loops-1} Let $(M,\tilde{g})$ be a pseudo-time orientable,
transverse, signature-type changing, $n$-dimensional ($n\geq2$)
manifold with a transverse radical, where $M_{L}=M\setminus(M_{R}\cup\mathcal{H})$
is globally hyperbolic. Assume that a Cauchy surface $S$ is a subset
of the neighborhood $U=\bigcup_{q\in\mathcal{H}}U(q)$ of $\mathcal{H}$,
i.e. $S\subseteq(U\cap M_{L})=\bigcup_{q\in\mathcal{H}}(U(q)\cap M_{L})$,
with $U(q)$ being constructed as in Theorem~\ref{Theo: Local Loops-1}.
Then for every point $p\in M$, there exists a pseudo-timelike loop
such that $p$ is a point of self-intersection.
\end{theorem}

\section{Pseudo-causal and pseudo-lightlike curves \label{sec:Pseudo-Timelike curves}}

In  $(M,g)$ be an $n$-dimensional manifold on which the metric signature changes from Riemannian,
 $(+,+,\text{\ensuremath{\cdots}},+)$, to Lorentzian,  $(\text{\textminus},+,\text{\ensuremath{\cdots}},+)$.
In such a setting, the conventional notion of timelike curves ceases to apply globally. From a point located in the Lorentzian region, the hypersurface of signature change may be reached within finite proper time. However, in the Riemannian region no meaningful notion of proper time exists, and hence curves there cannot be classified as timelike, spacelike, or null. Consequently, in manifolds with signature change
this gives rise to a novel notion of curves. In order to define those
curves we have to make a detour to draw upon the concept of the generalized
affine parameter.

\subsection{Properties of the generalized affine parameter\label{subsec:Properties of the GAP}}

In this section we introduce the notions of \textit{pseudo-timelike} and \textit{pseudo-spacelike} curves. To make these notions well-defined, we must distinguish genuine pseudo-timelike (respectively pseudo-spacelike) curves from those that asymptotically become lightlike as they approach
the hypersurface of signature change. For this purpose, the concept of a generalized affine parameter will play a central role. In particular, we require a suitable notion of completeness, ensuring that every $C^{1}$ curve of finite length with respect to such a parameter possesses an endpoint.
The use of generalized affine parameters in this context was first proposed by Ehresmann~\cite{Ehresmann}
and later developed further by Schmidt~\cite{Schmidt - A new definition of singular points in general relativity},
who employed them to characterize the completeness of arbitrary curves. The generalized affine parameter proves especially useful
quantity for probing singularities, because it can be defined for an
arbitrary curve, not necessarily a geodesic.

~

\begin{definition}[Generalized affine parameter]
\label{Definition-GAP-(Beem):Definition Generalized Affine Parameter}
Let $M$ be an $n$-dimensional manifold with an affine connection,
and let $\gamma\colon J\rightarrow M$ be a $C^{1}$ curve. A smooth vector field $V$ along $\gamma$ is a smooth map $V\colon J\rightarrow TM$
such that $V(t)\in T_{\gamma(t)}M$ for all $t\in J$. Such a 
vector field $V$ along $\gamma$ is called \emph{parallel} along
$\gamma$ if it satisfies the differential equation \[
\nabla_{\gamma'} V(t) = 0 \quad \text{for all } t \in J
\]
(see~\cite{Beem + Ehrlich + Easley - Global Lorentzian Geometry} for details).

\medskip
Fix $t_{0}\in J$ and  and set \(p_{0} = \gamma(t_{0})\). Choose a basis \(\{e_{1}, \dots, e_{n}\}\) of $T_{\gamma(t_{0})}M$. Let $E_{i}$ be the unique
parallel field along $\gamma$ with $E_{i}(t_{0})=e_{i}$ for $1\leq i\leq n$.
Then $\{E_{1}(t),E_{2}(t),\ldots,E_{n}(t)\}$ forms a basis for $T_{\gamma(t)}M$
for each $t\in J$. The tangent vector of $\gamma$ at $p_{0}$, can be expressed in this frame as a linear combination of the elements of the chosen basis with coefficients 
\[
V^{i}(t): \dot{\gamma}(t)=\sum_{i=1}^{n}\underset{V^{i}(t)E_{i}(t)}{\underbrace{V^{i}(t)E_{\gamma(t)i}}},
\]
where $V^{i}\colon J\longrightarrow\mathbb{R}$ for $1\leq i\leq n$.
Then the generalized affine parameter $\mu=\mu(\gamma,E_{1},\ldots,E_{n})$
of $\gamma(t)$ associated with this basis is given by 
\begin{equation}
\mu(t)=\int_{t_{0}}^{t}\sqrt{\mathop{\sum_{i=1}^{n}[V^{i}(t)]^{2}}}dt=\int_{t_{0}}^{t}\sqrt{\delta_{ij}V^{i}(t)V^{j}(t)}dt,\;t\in J.
\end{equation}
The assumption that \(\gamma\) is \(C^{1}\) ensures that the parallel fields \(\{E_{1}, \dots, E_{n}\}\) are well-defined. 
\end{definition}

~

Furthermore, we have\\

\begin{proposition}
~\textup{\cite{Hawking + Ellis - The large scale structure of spacetime}}
\label{prop:finite arc length in GAP}The curve $\gamma$ has a finite
arc-length in the generalized affine parameter $\mu=\mu(\gamma,E_{1},\ldots,E_{n})$
if and only if $\gamma$ has finite arc-length in any other generalized
affine parameter $\mu=\mu(\gamma,\bar{E}_{1},\ldots,\bar{E}_{n})$.
\end{proposition}

~

Note that the generalized affine parameter of a curve depends on the
chosen basis. Conceptually, one treats the parallel-transported basis vectors as if they were orthonormal with respect to a Riemannian metric, and then defines the ``length'' of \(\gamma(t)\) accordingly.  
In particular, if the metric \(g\) is positive definite, the generalized affine parameter associated with an orthonormal basis coincides with the usual arc-length.  

\noindent This notion of completeness succeeds in distinguishing precisely the behavior we intend to capture.  
A further advantage is that \(\mu\) can be defined for any \(C^{1}\) curve, including null curves as well as timelike or spacelike ones.  
Moreover, whenever a curve has unbounded proper length, its generalized affine parameter is necessarily unbounded as well~\cite{Beem + Ehrlich + Easley - Global Lorentzian Geometry}.

~

\begin{claim}
If the generalized affine parameter of a curve is finite for one choice of parallel frame, then it is finite for every such choice. Hence, for completeness it suffices to know whether the generalized affine parameter is finite or infinite, regardless of the chosen frame.
\end{claim}

\begin{proof}
Let $\{E_{i}\}$ and $\{\tilde{E}_{j}\}$ be two parallel bases of $T_{\gamma(t)}M$ along $\gamma$, and the components $V^{i}(t)$ with respect to another
basis are given by $\tilde{V}^{j}(t)=\sum_{i=1}^{n}A_{i}^{j}V^{i}(t)$. Then
\[
\dot{\gamma}(t)=\sum_{i=1}^{n}V^{i}(t)E_{i}(t)=\sum_{j=1}^{n}\tilde{V}^{j}(t)\tilde{E}_{j}(t),
\]
with $\tilde{V}^{j}(t)=\sum_{i=1}^{n}A_{i}^{j}V^{i}(t)$, where $A=(A_{i}^{j})$ is a constant invertible matrix.  
Hence also $V^{j}(t)=\sum_{i=1}^{n}a_{i}^{j}\tilde{V}^{i}(t)$ with $A^{-1}=(a_{i}^{j})$. Then the generalized affine parameters with respect to these basis are
$\mu(t)=\int_{t_{0}}^{t}\sqrt{\mathop{\sum_{i=1}^{n}[V^{i}(t)]}^{2}}dt$
and $\tilde{\mu}(t)=\int_{t_{0}}^{t}\sqrt{\mathop{\sum_{i=1}^{n}[\tilde{V}^{i}(t)]^{2}}}dt$.
By direct estimate,
\[
\left|\tilde{V}(t)\right|=\left|\sum_{i=1}^{n}A_{i}^{j}V^{i}(t)\right|\leq\sum_{i=1}^{n}\mid A_{i}^{j}\mid\mid V^{i}(t)\mid\leq\underset{ij}{\max}\sideset{\mid A_{i}^{j}\mid}{_{i=1}^{n}}\sum\mid V^{i}(t)\mid.
\]
\\ Applying Cauchy–Schwarz gives
\begin{multline*}
|\tilde{V}^{j}(t)|^{2} 
\leq \max_{i,j} |A_{i}^{j}|^{2}
   \underbrace{\left(\sum_{i=1}^{n} |V^{i}(t)| \right)^{2}}_{\left(\sum_{i=1}^{n} |V^{i}(t)| \cdot 1\right)^{2}} \\
\leq \max_{i,j} |A_{i}^{j}|^{2}
   \left(\sum_{i=1}^{n} |V^{i}(t)|^{2}\right)
   \left(\sum_{i=1}^{n} 1\right)
= n \cdot \max_{i,j} |A_{i}^{j}|^{2} \sum_{i=1}^{n} |V^{i}(t)|^{2}.
\end{multline*}

\noindent Summing over $j$ yields
\[
\sum_{j=1}^{n} |\tilde{V}^{j}(t)|^{2} 
\;\leq\; 
\sum_{j=1}^{n} \Biggl( n \cdot \max_{i,j} |A_{i}^{j}|^{2} \sum_{i=1}^{n} |V^{i}(t)|^{2} \Biggr) 
\;=\; 
n^{2} \cdot \max_{i,j} |A_{i}^{j}|^{2} \sum_{i=1}^{n} |V^{i}(t)|^{2}.
\]
On the other hand, we get 
$\sum_{j=1}^{n} |V^{j}(t)|^{2}
   \;\leq\;
   n^{2} \max_{i,j} |a_{i}^{j}|^{2}
   \sum_{i=1}^{n} |\tilde{V}^{i}(t)|^{2}.$ Combining both estimates yields
\[
\sum_{j=1}^{n} |\tilde{V}^{j}(t)|^{2} 
\leq n^{2} \cdot \max_{i,j} |A_{i}^{j}|^{2} \left( \sum_{i=1}^{n} |V^{i}(t)|^{2} \right)
\leq n^{4} \cdot \max_{i,j} |A_{i}^{j}|^{2} \cdot \max_{i,j} |a_{i}^{j}|^{2} \left( \sum_{i=1}^{n} |\tilde{V}^{i}(t)|^{2} \right).
\]
This is equivalent to
\[
\frac{1}{n^{2}\max_{i,j}|A_{i}^{j}|^{2}}
   \sum_{j=1}^{n}|\tilde{V}^{j}(t)|^{2}
   \;\leq\;
   \sum_{i=1}^{n}|V^{i}(t)|^{2}
   \;\leq\;
   n^{2}\max_{i,j}|a_{i}^{j}|^{2}
   \sum_{i=1}^{n}|\tilde{V}^{i}(t)|^{2},
\]
which implies
\[
\underbrace{\frac{1}{\sqrt{\,n^{2}\max_{i,j}|A_{i}^{j}|^{2}}}}_{c_{1}}
   \sqrt{\sum_{j=1}^{n}|\tilde{V}^{j}(t)|^{2}}
   \;\leq\;
   \sqrt{\sum_{i=1}^{n}|V^{i}(t)|^{2}}
   \;\leq\;
\underbrace{\sqrt{\,n^{2}\max_{i,j}|a_{i}^{j}|^{2}}}_{c_{2}}
   \sqrt{\sum_{i=1}^{n}|\tilde{V}^{i}(t)|^{2}}.
\]
Integrating, the generalized affine parameters
\[
\mu(t)=\int_{t_{0}}^{t}\sqrt{\sum_{i=1}^{n}|V^{i}(s)|^{2}}\,ds,
\quad
\tilde{\mu}(t)=\int_{t_{0}}^{t}\sqrt{\sum_{j=1}^{n}|\tilde{V}^{j}(s)|^{2}}\,ds
\]
satisfy
\begin{equation}{\label{eq: Estimate GAP}} 
c_{1}\tilde{\mu}(t)\le \mu(t)\le c_{2}\tilde{\mu}(t),
\end{equation}
for positive constants $c_{1},c_{2}$.  
Thus, $\mu(t)$ is finite if and only if $\tilde{\mu}(t)$ is finite.

\end{proof}

\subsection{Application of the generalized affine parameter in a signature-type
changing manifold\label{subsec:Application-of-the GAP}}

Let $M=M_{L}\cup\mathcal{H}\cup M_{R}$ be an $n$-dimensional transverse
type-changing singular semi-Riemannian manifold with a type-changing
metric $g$, and $\mathcal{H}:=\{q\in M\!\!:g\!\!\mid_{q} is\;degenerate\}$
the locus of signature change. We further assume that one component,
$M_{L}$, of $M\setminus\mathcal{H}$ is Lorentzian and the other
one, $M_{R}$, is Riemannian.

~

\begin{definition}[Pseudo-lightlike curve]
Given a continuous and differentiable curve $\gamma\colon[a,b]\longrightarrow M$,
with $[a,b]\subset\mathbb{R}$, where $-\infty<a<b<\infty$. Then
the curve $\gamma=\gamma^{\mu}(u)=x^{\mu}(u)$ is a pseudo-lightlike
curve if
~

\begin{itemize}
\item its tangent vector field in the Lorentzian component $M_{L}$ is null,
\item its tangent vector field in the Riemannian component $M_{R}$ is arbitrary.
\end{itemize}
\end{definition}

~

\begin{figure}[h]
\centering
\includegraphics[width=0.6\textwidth]{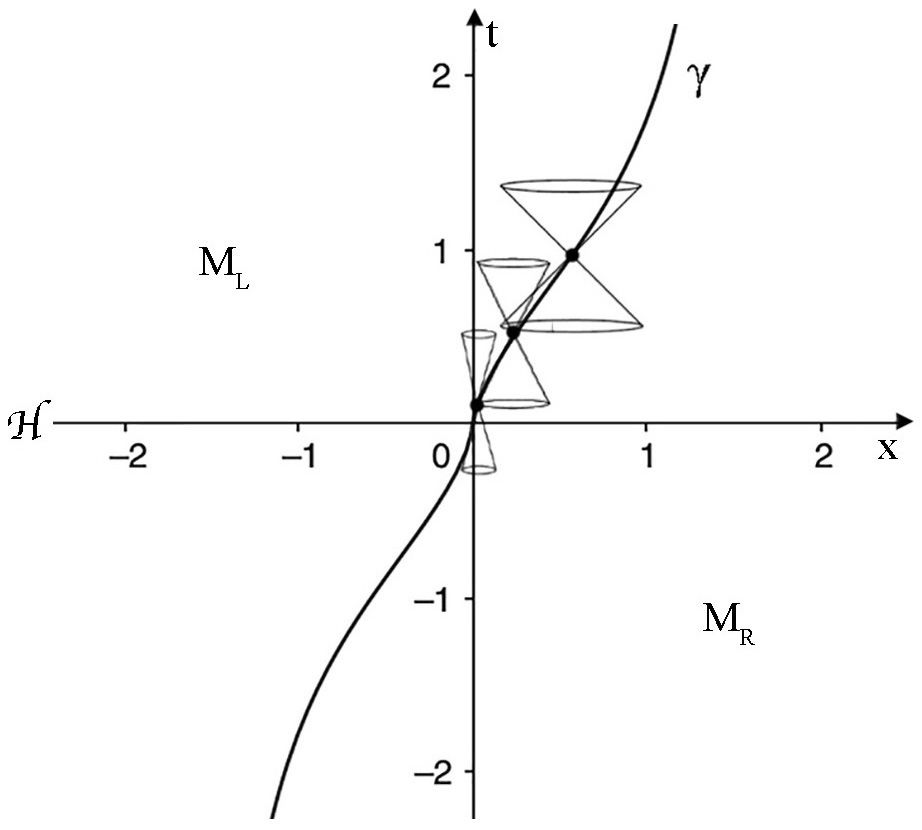}
\caption{The curve $\gamma$ is not pseudo-timelike since it approaches a null vector at the locus of signature change. This curve is asymptotically lightlike.}
\label{fig:asymptotically lightlike curve}
\end{figure}

A similar definition applies for \textit{\emph{a}}\emph{ }\textbf{\textit{\emph{pseudo-causal
curve}}}. Note that an analogous definition for pseudo-timelike and
pseudo-spacelike curves turns out to be problematic as the definition
would also include curves that asymptotically become lightlike as
they approach $\mathcal{H}$, see Figure~\ref{fig:asymptotically lightlike curve}.

~ 

\begin{example}
\label{Example Generalized Affine Parameter} Consider the metric $g = t\, (dt)^2 + (dx)^2$ on $\mathbb{R}^{2}$, and the non-parametrized, non-geodesic curve $\gamma$ defined by 
\[
\tan x = \frac{2}{3} \, \mathrm{sgn}(t)\, |t|^{3/2}, \quad -\frac{\pi}{2} < x < \frac{\pi}{2}.
\] 
We can rearrange this equation to isolate $x$ (see Figure~\ref{fig:GAP-Curve}):

~

\[
\frac{3}{2} \tan x = \operatorname{sgn}(x) \, \biggl|\frac{3}{2} \tan x\biggr| = \operatorname{sgn}(t) \, |t|^{\frac{3}{2}}
\]

\[
\Longleftrightarrow 
\underbrace{\operatorname{sgn}(t)\, |t|}_{t} = \operatorname{sgn}(x) \, \Bigl|\frac{3}{2} \tan x\Bigr|^{\frac{2}{3}}=t.
\]

\textbf{~}\\ Next,\footnote{Only the coordinate transformation for $T$ is taken from~\cite{Dray - General Relativity and Signature Change}. The example itself, including the choice of the curve $\gamma$ and the analysis of $\gamma$, are our original contributions.} we reintroduce the coordinate transformation suggested in~\cite{Dray - General Relativity and Signature Change}, expressed as
\[
T=\int_{0}^{t}\sqrt{\mid\tilde{t}\mid}d\tilde{t}=\frac{2}{3}\sqrt{\mid t\mid}^{3}\cdot\textrm{sgn}(t).\]

\textbf{~}\\ This yields the metric $g=\textrm{sgn}(T)(dT)^{2}+(dx)^{2}$,
and for the curve $\gamma$ we have $T=\tan x$. Thus, in
the $(T,x)$-coordinates, $\gamma$ is simply the $\tan$-function, with
derivative $dT/dx = 1/\cos^2 x$. Consequently, $\gamma$ is
timelike in $M_{L}$ , approaching the light cone from timelike infinity
and touching it tangentially as $T\rightarrow0$ (where
the derivative becomes $\frac{1}{\cos^{2}(0)}=1$ in the limit). These are precisely the types of curves we aim to exclude in our definition. Note that unlike $(T,x)$, the $(t,x)$-coordinates cover the entire manifold $M$.

\begin{figure}[h]
\centering
\includegraphics[width=0.75\textwidth]{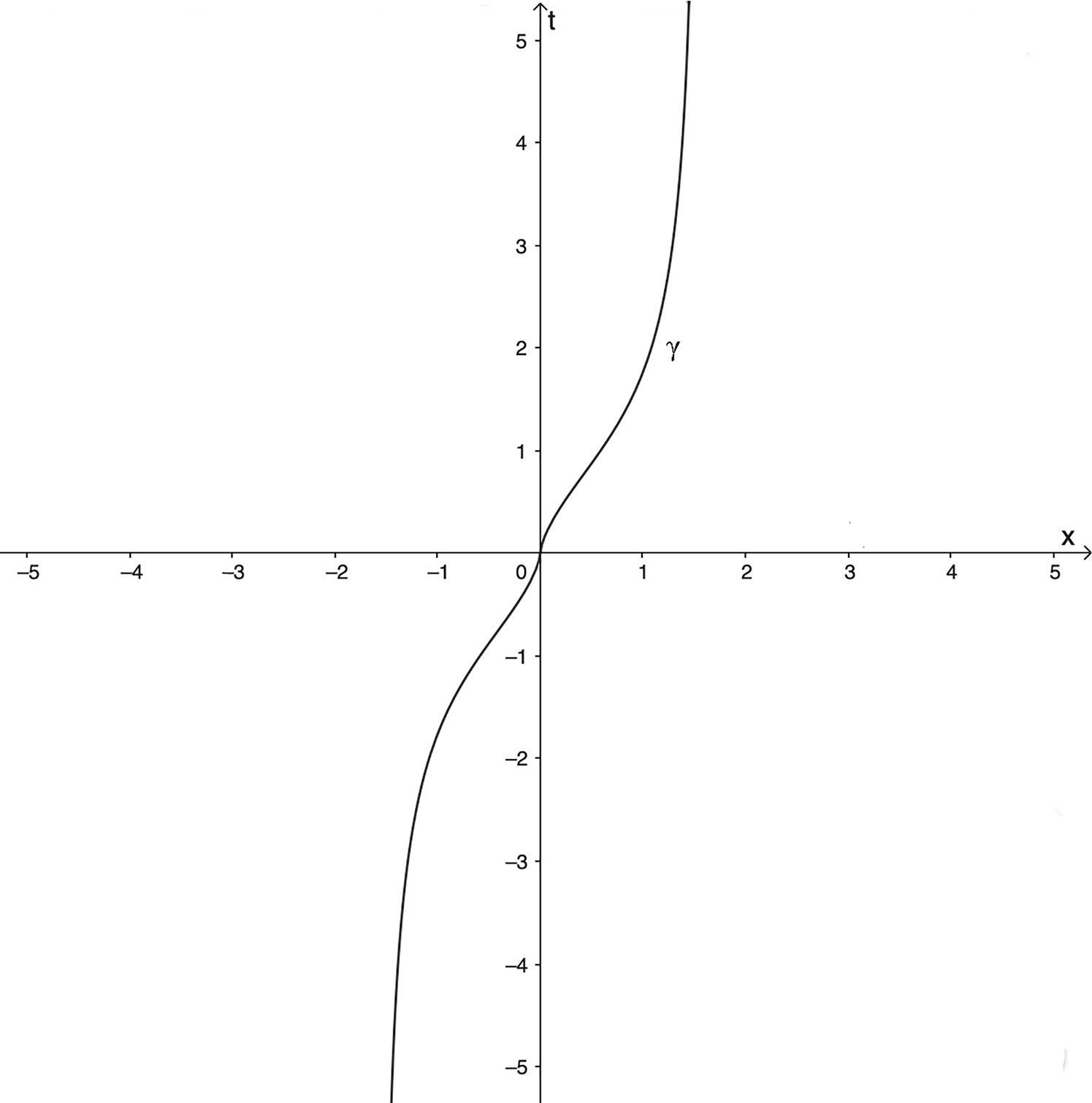}
\caption{The curve defined by $t=\textrm{sgn}(x)\cdot\left(\mid\frac{3}{2}\tan x\mid\right)^{\frac{2}{3}}$.}
\label{fig:GAP-Curve}
\end{figure}

Moreover, if the curve $\gamma=(T(s),x(s))$ is parametrized by arc
length $s$, then in the $(t,x)$-coordinate system both $\frac{dx}{ds}$
and $\frac{dt}{ds}$ diverge in $M_{L}$, since $-1=-(\frac{dT}{ds})^{2}+(\frac{dx}{ds})^{2}$:

~

\[
-1 =  -\left(\frac{dT}{dx} \frac{dx}{ds}\right)^2 +\left(\frac{dx}{ds}\right)^2 = \left(-\left(\frac{d\tan x}{dx}\right)^2 + 1\right) \left(\frac{dx}{ds}\right)^2 = \left(-\frac{1}{\cos^4 x} + 1\right) \left(\frac{dx}{ds}\right)^2,
\]
which is equivalent to $(\frac{dx}{ds})^{2}=\frac{-1}{(-\frac{1}{\cos^{4}x}+1)}$. And this implies

\[
\underset{x\rightarrow0}{\lim}\frac{dx}{ds}=\underset{x\rightarrow0}{\lim}\pm\sqrt{\frac{-1}{(-\frac{1}{\cos^{4}x}+1)}}=\pm\infty.
\]
Similarly, we get
\begin{align}
-1 &= \left(-1 + \frac{1}{\left(\frac{dT}{dx}\right)^{2}}\right) \left(\frac{dT}{ds}\right)^{2} \notag\\
&= \left(-1 + \underset{\frac{1}{(\frac{d\tan x}{dx})^{2}}}{\underbrace{\cos^{4}x}}\right) \left(\frac{dT}{dt}\right)^{2} \left(\frac{dt}{ds}\right)^{2} \notag\\
&= \left(-1 + \frac{1}{(1 + \tan^{2}x)^{2}}\right) |t| \left(\frac{dt}{ds}\right)^{2} \notag\\
&= \left(-1 + \frac{1}{(1 + T^2)^2}\right) |t| \left(\frac{dt}{ds}\right)^2 \notag\\
&= \left(-1 + \frac{1}{(1 + \frac{4}{9} |t|^3)^2}\right) |t| \left(\frac{dt}{ds}\right)^2.
\end{align}

\noindent This is equivalent to
\begin{equation}
\left(\frac{dt}{ds}\right)^{2} = \frac{-1}{\left(-1 + \frac{1}{\left(1 + \frac{4}{9} |t|^3 \right)^2}\right) |t|},
\end{equation}
and then yields
\begin{equation}
\lim_{t \to 0} \frac{dt}{ds} = \lim_{t \to 0} \pm \sqrt{\frac{-1}{\left(-1 + \frac{1}{(1 + \frac{4}{9} |t|^3)^2}\right) |t|}} = \pm \infty.
\end{equation}

\end{example}

~

While the components of $\gamma'$ do not diverge in the $(T,x)$-coordinate
system, both $\frac{dx}{ds}$ and $\frac{dt}{ds}$ diverge in $M_{L}$
in the $(t,x)$-coordinates. Because of this dependency of coordinates
the criterion of divergence is not useful for defining pseudo-timelike
and pseudo-spacelike curves. That is where the coordinate-independent
generalized affine parameter comes into play.

~

\begin{definition}[Pseudo-timelike  
curve]
\textbf{\label{Definition Pseudo-timelike curve}} Let $M=M_{L}\cup\mathcal{H}\cup M_{R}$ be an $n$-dimensional
transverse type-changing singular semi-Riemannian manifold, $g$ 
a type-changing metric, and $\mathcal{H}:=\{q\in M\!\!:g\!\!\mid_{q} is\;degenerate\}$
the locus of signature change. We assume that one component,
$M_{L}$, of $M\setminus\mathcal{H}$ is Lorentzian, and the other
component, $M_{R}$, is Riemannian.\\ Let $\gamma\colon [a,b] \to M$ be a continuous and differentiable
curve, with $[a,b]\subset\mathbb{R}$, $-\infty<a<b<\infty$. We call $\gamma=\gamma^{\mu}(u)=x^{\mu}(u)$
in $M$ a \textit{\emph{pseudo-timelike}}\emph{ (respectively,
pseudo-spacelike)} curve if

1. $\operatorname{Im}(\gamma) \cap M_{L} \neq \varnothing$, i.e. $\gamma$ has image points in the Lorentzian region, and

2. for every generalized affine parametrization of $\gamma$ in $M_{L}$ there exists $\varepsilon>0$ such that  
\[
g(\gamma',\gamma') < -\varepsilon 
\quad \text{(respectively, } g(\gamma',\gamma') > \varepsilon \text{)}.\footnote{Since Definition~\ref{Definition-GAP-(Beem):Definition Generalized Affine Parameter}
is already independent of a choice of coordinates and instead refers
to a (generally anholonomic) basis, the above Definition~\ref{Definition Pseudo-timelike curve}
is also coordinate independent. The independence of Definition~\ref{Definition Pseudo-timelike curve}
from the choice of this basis is a direct consequence of Proposition~\ref{prop:finite arc length in GAP}.
In particular, in the case of a basis change we just relegate to the
estimate~(\ref{eq: Estimate GAP}).}
\]
\end{definition}

~

\begin{example}
\textit{\label{Example-1}} Revisiting Example~\ref{Example Generalized Affine Parameter},
we find that both coordinate vector fields, $\frac{\partial}{\partial T}$
and $\frac{\partial}{\partial x}$, are covariantly constant in $M_{L}$
and $M_{R}$ (this is because the Christoffel symbols all vanish in
the $(T,x)$-coordinate system). Hence, we can parallel transport
$\frac{\partial}{\partial T}$ and $\frac{\partial}{\partial x}$
along any curve in $M_{L}$ and $M_{R}$, with the transport being
path-independent (no anholonomy).\\ Since we aim at parametrizing
the curve $\gamma$ by the generalized affine parameter $\mu$ with
respect to the coordinate vector fields $\frac{\partial}{\partial T}=\frac{1}{\sqrt{\mid t\mid}}\frac{\partial}{\partial t}$
and $\frac{\partial}{\partial x}$ we are able to start with an arbitrary
parametrization. Hence, let $\gamma(t)=(T(t),x(t))$ be parametrized
by $t$, and then $\dot{\gamma}(t)=\frac{dT}{dt}\frac{\partial}{\partial T}+\frac{dx}{dt}\frac{\partial}{\partial x}$.
By means of Definition~\ref{Definition-GAP-(Beem):Definition Generalized Affine Parameter}
we immediately get 
\[
V^{0}(t)=\frac{dT}{dt}=\sqrt{\mid t\mid}
\]
 and
\[
V^{1}(t)=\frac{dx}{dt}=\frac{d}{dt}\arctan\left(\frac{2}{3}\sqrt{\mid t\mid}^{3}\textrm{sgn}(t)\right).
\]
In $M_{L}$ this yields 
\[
V^{1}(t)=\frac{\sqrt{\mid t\mid}}{1+\frac{4}{9}\mid t\mid^{3}}.
\]
Consider now $\tilde{\gamma}(t(s))=\gamma(s)$, in which $\tilde{\gamma}$
is related to the curve $\gamma$ by reparametrization of $\gamma$
by $t$. With this notation we have the basis fields $E_{\tilde{\gamma}(t),0}=\frac{1}{\sqrt{\mid t\mid}}\frac{\partial}{\partial t}$
and $E_{\tilde{\gamma}(t),1}=\frac{\partial}{\partial x}$ along $\tilde{\gamma}$.
The reparametrized curve $\tilde{\gamma}(t(s))$ also gives 
\[
\dot{\tilde{\gamma}}(t)=V^{i}(t)E_{\tilde{\gamma}(t),i}=\frac{\partial}{\partial t}+\frac{dx}{dt}\frac{\partial}{\partial x}.
\]
\textbf{~}\\ The Definition~\ref{Definition-GAP-(Beem):Definition Generalized Affine Parameter}
for the generalized affine parameter gives

\[
\frac{d\mu}{dt}=\sqrt{(V^{0}(t))^{2}+(V^{1}(t))^{2}}=\sqrt{\mid t\mid+\frac{\mid t\mid}{(1+\frac{4}{9}\mid t\mid^{3})^{2}}}.
\]
It now follows easily that for the reparametrization of $\hat{\gamma}(t)$
by the generalized affine parameter $\mu$ (i.e. $\hat{\gamma}(\mu(t))=\tilde{\gamma}(t)$)
we have in $M_{L}$:
\[
g(\dot{\hat{\gamma}}(\mu(t)),\dot{\hat{\gamma}}(\mu(t)))=g(\frac{d\hat{\gamma}(\mu(t))}{d\mu},\frac{d\hat{\gamma}(\mu(t))}{d\mu})=g(\frac{1}{\frac{d\mu}{dt}}\dot{\hat{\gamma}}(t),\frac{1}{\frac{d\mu}{dt}}\dot{\hat{\gamma}}(t))
\]
\[
=\frac{g(\frac{\partial}{\partial t}+\frac{dx}{dt}\frac{\partial}{\partial x},\frac{\partial}{\partial t}+\frac{dx}{dt}\frac{\partial}{\partial x})}{(\frac{d\mu}{dt})^{2}}=\frac{t+(\frac{dx}{dt})^{2}}{\mid t\mid+\frac{\mid t\mid}{(1+\frac{4}{9}\mid t\mid^{3})^{2}}}=\frac{t+\frac{\mid t\mid}{(1+\frac{4}{9}\mid t\mid^{3})^{2}}}{\mid t\mid+\frac{\mid t\mid}{(1+\frac{4}{9}\mid t\mid^{3})^{2}}}.
\]
Taking the limit 
\[
\underset{t\rightarrow0^{-}}{\lim}\frac{t+\frac{\mid t\mid}{(1+\frac{4}{9}\mid t\mid^{3})^{2}}}{\mid t\mid+\frac{\mid t\mid}{(1+\frac{4}{9}\mid t\mid^{3})^{2}}}=0
\]
reveals that the curve $\gamma$ is not pseudo-timelike as it does
not meet the $\varepsilon$-requirement of Definition~\ref{Definition Pseudo-timelike curve}.
\end{example}

~

\textbf{~}\\ In Section~\ref{subsec:Properties of the GAP} we repeatedly
rather vaguely referred to the concept of a timelike (or spacelike,
respectively) curve that asymptotically becomes lightlike. The above
example highlights how the notion of ``asymptotically lightlike''
should be understood. A timelike (or spacelike, respectively) curve
in $M_{L}$ that is not pseudo-timelike (or pseudo-spacelike, respectively)
can be thus specified as asymptotically lightlike.

~

\begin{example}
\textit{\label{Example-2.}}\textbf{ }Finally, if we modify the previously
discussed curve $\gamma$ by keeping the $t$-coordinate but stating
$x=0$, we get a curve $\alpha$. With the same notation as above,
we then get $V^{0}(t)=\sqrt{\mid t\mid}$, $V^{1}(t)=0$ and $\frac{d\mu}{dt}=\sqrt{\mid t\mid}$.
Hence, this results in $g(\dot{\hat{\alpha}}(\mu(t)),\dot{\hat{\alpha}}(\mu(t)))=\frac{1}{\mid t\mid}g(\frac{\partial}{\partial t},\frac{\partial}{\partial t})=\frac{t}{\mid t\mid}=-1$
in the Lorentzian region $M_{L}$. The curve $\alpha$ is pseudo-timelike
as it obviously does meet the $\varepsilon$-requirement of Definition~\ref{Definition Pseudo-timelike curve}.
\end{example}

~

This disquisition makes it clear why the notion of the generalized affine parameter is necessary and useful in order to define pseudo-timelike and pseudo-spacelike curves. If we were to loosen the requirement in Definition~\ref{Definition Pseudo-timelike curve} by replacing ``\textit{for every generalized affine parametrization of $\gamma$ in $M_{L}$}'' with ``\textit{for every affine parametrization of $\gamma$ in $M_{L}$}'', then no curve that is timelike in the Lorentz sector and reaches the hypersurface $\mathcal{H}$ would be pseudo-timelike throughout the entire manifold $M$. (However, this statement applies only to curves that actually reach the hypersurface $\mathcal{H}$. Timelike curves that lie entirely within $M_{L}$ and maintain a ``distance'' from $\mathcal{H}$ due to a tubular neighborhood within $M_{L}$ also satisfy the relaxed condition, as they only have affine parameters with $g(\gamma', \gamma') = \text{const} < 0$.)

~

Similarly, any timelike curve in $M_{L}$ would meet the requirements
of a pseudo-timelike curve if we modified the definition by requesting
``\textit{for a suitable parametrization of $\gamma$ in $M_{L}$}''
instead of ``\textit{for every generalized affine parametrization
of $\gamma$ in $M_{L}$}''. In this regard, the concept of the generalized
affine parameter is the right tool to tell apart suitable from unsuitable
curves for the definition of pseudo-timelike and pseudo-spacelike
curves. 

~

Interestingly, our rationale for the new definition of a pseudo-timelike
curve is reminiscent of the analysis undertaken in~\cite{Ling}.
In Section 2 of~\cite{Ling} the distinction between \textit{causal
curves}, \textit{timelike almost everywhere curves} and \textit{timelike
curves} is introduced in which the latter one is defined as follows:
A timelike curve is a causal curve $\gamma\colon I\longrightarrow M$
such that $g(\gamma',\gamma')<-\varepsilon$ almost everywhere for
some $\varepsilon>0$.\textbf{~}\\ The author illustrates the situation
in his Figure 1 which contrasts a timelike curve with a timelike almost
everywhere curve. The latter one can not be viewed as a timelike curve
because it approaches a null vector at its break point. Compared to
our setting, however, the culprit here is that the curve is not differentiable
at the breaking point. However, if we were to make the curve differentiable
by bending its upper section, it still wouldn't be timelike. Its restriction
to the lower region before the inflection point is timelike, but it
cannot be extended upwards into a timelike curve. Now, imagine we
are not in Minkowski spacetime, but instead, a signature-type change
occurs at the (former) inflection point so that the ``upper half'' of
the space becomes Euclidean (in this case, the figure would correspond
to the $(T,x)$-coordinates, not the $(t,x)$-coordinates, in the
toy model). In this scenario, the curve restricted to the Lorentz
sector would be timelike, but after the signature change is ``reversed'',
it cannot be extended upwards into a timelike curve. In this sense,
the entire curve in the signature-changing version of this example
is not pseudo-timelike.

~ 
\begin{remark}
The introduction of pseudo-timelike curves is not merely a technical tool. If a timelike trajectory were to become asymptotically null near the degenerate hypersurface, this would correspond to a massive particle being accelerated to arbitrarily close to the speed of light. Such behaviour would require a divergent energy input and is therefore unphysical. Pseudo-timelike curves exclude this scenario and ensure that the worldlines of massive particles remain genuinely timelike.
\end{remark}

~

Now we can slightly modify the definition of a (simply) closed curve
in order for it to correctly apply to signature-type changing singular
semi-Riemannian manifolds $M$ with a metric $g$:

~

\begin{definition}[Chronology-violating
curve]
\label{Definition Chronology violating curve} A smooth, \textit{pseudo-timelike curve} $\gamma\colon I\longrightarrow M$
is said to be chronology-violating if there exists a subset of
 $\gamma[I]$ homeomorphic to $S^{1}$ such that there are at least two parameters
$s_{1},s_{2}\in I$ with $\gamma(s_{1})=\gamma(s_{2})$. Moreover,
$\gamma$ must belong to one of the following two classes:\footnote{Note that this means that there must be at least one such subset to
fulfill this definition.}\\

\begin{enumerate}
\item The pseudo-timelike curve $\gamma$ is periodic, i.e the image $\gamma[I]$
is homeomorphic to $S^{1}$. Moreover, for $s_{1}$, $s_{2}\in I$
the associated tangent vectors, $\gamma'(s_{1})$ and $\gamma'(s_{2})$,
are timelike and positively proportional. We denote this type of curve
as \textit{closed pseudo-timelike curve.}
\item The curve $\gamma$ intersects itself for $s_{1}$, $s_{2}\in I$
and the associated tangent vectors, $\gamma'(s_{1})$ and $\gamma'(s_{2})$,
are timelike whereas the tangent directions are not necessarily the
same (i.e. they do not need to be positively proportional). This type
of curve is said to contain a \textit{loop}. 
\end{enumerate}
\end{definition}

\section{Global structure of signature-type changing semi-Riemannian manifolds}

First, let us revisit the definitions of the following concepts related
to manifold orientability.

~

\begin{definition}
~\cite{Bott + Tu}  A smooth $n$-dimensional manifold $M$ is\textit{
orientable} if and only if it has a smooth global nowhere vanishing
$n$-form (also called a top-ranked form).\footnote{An orientation of $M$ is the choice of a continuous pointwise orientation, i.e. the specific choice of a global nowhere vanishing $n$-form.}
\end{definition}

~

For a differentiable manifold to be orientable all that counts is
that it admits a global top-ranked form - it is not important which
specific top-ranked form is selected. 

~ 

To ensure thoroughness, we also want to mention the definition of
parallelizability, which likewise does not involve any metric and
is therefore again applicable to manifolds with changing signature
types. It is well-known that a manifold $M$ of dimension $n$ is
defined to be parallelizable if there are $n$ global vector fields
that are linearly independent at each point. We define it similarly
to the approach in~\cite{Bredon}:

~

\begin{definition}
A smooth $n$-dimensional manifold $M$ is \textit{parallelizable}
if there exists a set of smooth vector fields $\{V,E_{1},\ldots,E_{n-1}\}$
on $M$, such that at every point $p\in M$ the tangent vectors $\{V(p),E_{1}(p),\ldots,E_{n-1}(p)\}$
provide a basis of the tangent space $T_{p}M$. A specific choice
of such a basis of vector fields on $M$ is called an \textit{absolute
parallelism} of $M$.
\end{definition}

~

Equivalently, a manifold $M$ of dimension $n$ is parallelizable
if its tangent bundle $TM$ is a trivial bundle, so that the associated
principal bundle of linear frames has a global section on $M$, i.e.
the tangent bundle is then globally of the form $TM\simeq M\times\mathbb{R}^{n}$.
Moreover, it is worth pointing out that orientability and also parallelizability
are \textit{differential topological properties} which do not depend
on the metric structure, but only on the topological manifold with
a globally defined differential structure.

~ 

\begin{remark}
It is worth mentioning that given an absolute parallelism of $M$,
one can use these $n$ vector fields to define a basis of the tangent
space at each point of $M$ and thus one can always get a frame-dependent
metric $g$ by defining the frame to be orthonormal. Moreover, the
special orthogonal group, denoted $SO(n,\mathbb{R})$, acts naturally
on each tangent space via a change of basis, it is then possible to
obtain the set of all orthonormal frames for $M$ at each point qua
the oriented orthonormal frame bundle of $M$, denoted $F_{SO}(M)$,
associated to the tangent bundle of $M$. 
\end{remark}

~

The next three definitions, however, depend not only on the underlying
manifold but also on its specific type-changing metric $g$. For our
purpose, let $(M,g)$ be a smooth, signature-type changing manifold
(possibly with boundary).

~

\begin{definition}[Pseudo-timelike]
\label{def:Pseudo-timelike} A vector field $V$ on a signature-type changing manifold $(M,g)$ is \textit{pseudo-timelike} if and only if its integral curves are pseudo-timelike (hence, in particular, $V$ is timelike in $M_{L}$).\footnote{Keep in mind that a timelike vector field is a vector field $V$ on
a spacetime manifold $(M,g)$ where the vectors at every point are
timelike, meaning $g(V(p),V(p))<0$ for all points $p$ on the manifold.}
\end{definition}

~

\begin{definition}[Pseudo-time orientable]
\label{def:Pseudo-time-orientable} 
A signature-type changing manifold $(M,g)$ is \textit{pseudo-time
orientable} if and only if the Lorentzian region $M_{L}$ is time
orientable.\footnote{A pseudo-time orientation of such a manifold $(M,g)$ corresponds
to the specific choice of a continuous non-vanishing pseudo-timelike
vector field $V$ on $M$.}
\end{definition}

~

\begin{lemma}
A singular semi-Riemannian manifold \textup{$(M,g)$} is pseudo-time orientable if and only if there exists a vector field $X \in \mathfrak{X}(M)$ that is pseudo-timelike.
\end{lemma}

\begin{proof}

$"\Longrightarrow"$ Let a singular semi-Riemannian manifold $(M,g)$
be pseudo-time orientable. This means the Lorentzian region $M_{L}$
is time orientable. A Lorentzian manifold is time-orientable if there
exists a continuous timelike vector field. Accordingly, there must
exist a continuous timelike vector field $X \in \mathfrak{X}(M_{L})$
in the Lorentzian region. As per Definition~\ref{def:Pseudo-timelike},
a vector field $X$ in a signature-type changing manifold is pseudo-timelike
if and only if $X$ is timelike in $M_{L}$ and its integral curve
is pseudo-timelike; this means that $X$ is allowed to vanish on $M_{R}$.\footnote{In this part of the proof, the only thing that matters is whether
the ``pseudo-timelike vector field'' is
allowed to vanish on $M_{R}$. This question is independent of whether
the ``generalized affine parameter'' condition
is required in $M_{L}$, because the issue of whether the vector field
``is allowed to vanish on $M_{R}$'' concerns only its
``magnitude'', while the ``generalized
affine parameter'' condition pertains solely to its ``direction''
(specifically, that the vector field is not asymptotically lightlike).} Hence, we can extend the vector field $X$ arbitrarily to all of
$M$, and per definition $X \in \mathfrak{X}(M)$ is pseudo-timelike.

~

$"\Longleftarrow"$ Let $X\in\mathfrak{X}(M)$ be a pseudo-timelike
vector field in a singular semi-Riemannian manifold $(M,g)$. Hence,
as per Definition~\ref{def:Pseudo-timelike}, $X$ is timelike in
$M_{L}$. A Lorentzian manifold is time-orientable if and only if
there exists a timelike vector field. Since $X$ is a timelike vector
field on $M_{L}$, the Lorentzian region $M_{L}$ is time-orientable.
Then, according to Definition~\ref{def:Pseudo-time-orientable},
the signature-type changing manifold $(M,g)$ is pseudo-time orientable.\textit{ }
\end{proof}

According to that, such a definition of a pseudo-time orientation
is possible if $M_{L}$ admits a globally consistent sense of time,
i.e. if in $M_{L}$ we can continuously define a division of non-spacelike
vectors into two classes. For a transverse, signature-type changing
manifold (with a transverse radical), this definition arises naturally
because, in $M_{R}$, all vectors can be considered spacelike. Additionally,
all non-spacelike vectors on $\mathcal{H}$ are lightlike.\footnote{Note that this applies generally, including in the case of a tangent
radical, since there are no timelike vectors on $\mathcal{H}$. However,
the subsequent division into two classes requires a transverse radical.} In the case that $Rad_{q}\cap T_{q}\mathcal{H}=\{0\}$ $\forall q\in\mathcal{H}$,
these lightlike vectors can be naturally divided into two classes:
those pointing towards $M_{L}$ and those pointing towards $M_{R}$.

~

\begin{example}

\begin{figure}[h]
\centering
\includegraphics[width=0.5\textwidth]{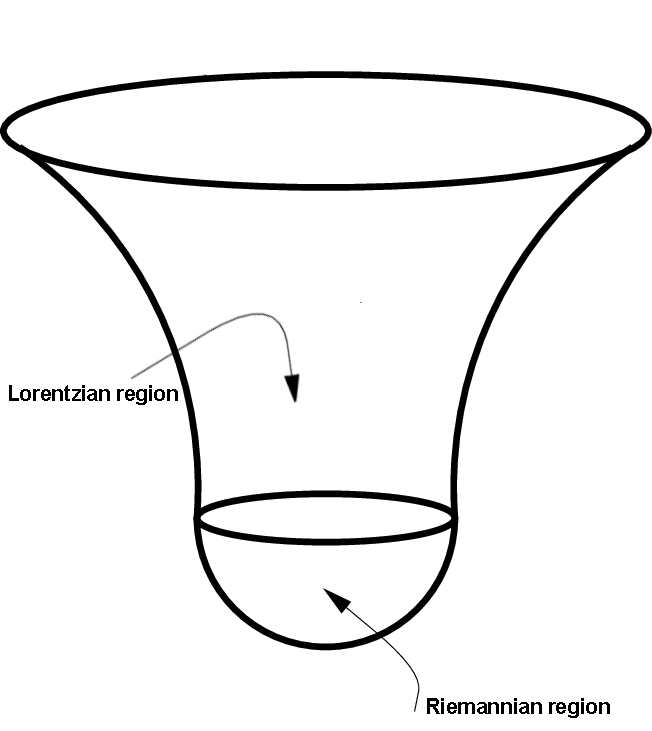}
\caption{Riemannian
and Lorentzian region in the Hartle-Hawking no-boundary model.}
\label{fig:No-Boundary Modell}
\end{figure}

\label{Example-7}Consider the classic type of a spacetime $M$ with
signature-type change which is obtained by cutting an $S^{4}$ along
its equator and joining it to the corresponding half of a de Sitter
space, Figure~\ref{fig:No-Boundary Modell}. The de Sitter spacetime
is time-orientable~\cite{Rodrigues - de Sitter}, hence $M$ is pseudo-time
orientable.
\end{example}

~

\begin{definition}[Pseudo-space orientable]
\label{def:Pseudo-space-orientable} 
A signature-type changing manifold $(M,g)$ of dimension $n$ is \textit{pseudo-space
orientable} if and only if it admits a continuous non-vanishing spacelike
$(n-1)$-frame field on $M_{L}$. This is a set of $n-1$ pointwise
orthonormal spacelike vector fields on $M_{L}$.\footnote{A pseudo-space orientation of a manifold $(M,g)$ corresponds to the
specific choice of a continuous non-vanishing field of orthonormal
spacelike $(n-1)$-beins on $M_{L}$.}
\end{definition}

~

\begin{proposition}
\label{Proposition. parallelizable-orientable}\textup{\cite{Lee - Introduction to Smooth Manifolds}}
Every parallelizable manifold $M$ is orientable.
\end{proposition}

~

In Lorentzian geometry the fact of $M$ being time-orientable and
space-orientable implies that $M$ is orientable~\cite{Hall}. The
proposition below illustrates that this result from Lorentzian geometry
cannot be applied to signature-type changing manifolds.

~

\begin{proposition}
\label{prop:orientable+timeorientable-not_orientable} Even if a transverse,
signature-type changing manifold $(M,g)$ with a transverse radical
is pseudo-time orientable and pseudo-space orientable, it is not necessarily
orientable.
\end{proposition}

\begin{proof}
Consider an arbitrary manifold of $\dim(M)=2$ with a change of signature,
for which the conditions of Proposition~\ref{prop:orientable+timeorientable-not_orientable}
are given (in higher dimensions, the same idea can be carried out
through a trivial augmentation of dimensions). In case this manifold
is non-orientable, there is nothing to show. However, if it is orientable,
cut out a disk from the Riemannian sector and replace it with a crosscap,
equipped with any Riemannian metric. In a tubular neighborhood of
the cutting line, construct a Riemannian metric that mediates between
the metrics of the crosscap and the rest (this is possible due to
the convexity of the space formed by all Riemannian metrics). This
surgical intervention results in the transition to a non-orientable
manifold with a change of signature. Since the intervention is limited
to the Riemannian sector, the conditions of the proposition remain
unaffected. Thus Proposition~\ref{prop:orientable+timeorientable-not_orientable}
is proven.
\end{proof}

~ 

\begin{remark}
One can always ``switch'' between non-orientability and orientability using the crosscap. Starting with an orientable manifold, one transitions to non-orientable by replacing a crosscap (if already present) with a disk. If no crosscap is present, such a transition occurs by replacing a disk with a crosscap.
\end{remark}

~

\begin{example}\label{Example-9}
The M\"{o}bius strip $\mathbb{M}$ has a non-trivial vector bundle structure over $S^{1}$, which means that the bundle cannot be trivialized globally. Specifically, the M\"{o}bius strip is a line bundle over $S^{1}$ with a non-trivial twist.\footnote{The M\"{o}bius strip is particularly interesting because it can be found on any arbitrary non-orientable surface. Additionally, any Lorentzian manifold $\mathbb{M} \times \mathbb{R}^{n}$ based on the M\"{o}bius strip crossed with $\mathbb{R}^{n}$ either fails to be time-orientable or space-orientable~\cite{Geroch - Global structure of spacetimes}.} Hence, $\mathbb{M}$ is neither parallelizable nor orientable.

\textbf{~}\\To see this, consider the M\"{o}bius strip $\mathbb{M} = \mathbb{R} \times \mathbb{R}/\sim$ with the identification $(t,x) \sim (\tilde{t},\tilde{x}) \Longleftrightarrow (\tilde{t}, \tilde{x}) = ((-1)^{k}t, x + k), k \in \mathbb{Z}$. Notice that the identification has no bearing on proper subsets of $((-1)^{k}t, x + k), k \in \mathbb{Z}$, and the fibre $\mathbb{R}$ is a vector space.

\textbf{~}\\As $\mathbb{M}$ is a fiber bundle over the base space $S^{1}$, a section of that fiber bundle must be a continuous map $\sigma: S^{1} \longrightarrow \mathbb{M}$ such that $\sigma(x) = (h(x), x) \in \mathbb{M}$. For $\sigma$ to be continuous, $h$ must satisfy $-h(0) = h(k)$. The intermediate value theorem guarantees that there is some $\tilde{x} \in [0, k]$ such that $h(\tilde{x}) = 0$. This means that every section of $\mathbb{M}$ intersects the zero section, and the sections that form a basis for the fibre are not non-zero everywhere.
\end{example}

~

\begin{definition}
A \textbf{pseudo-spacetime} is a $4$-dimensional, pseudo-time oriented,
semi-Riemannian manifold with a type-changing metric. 
\end{definition}

~

\begin{proposition}
Let $(\mathbb{R}^{n}, g)$ be a transverse, signature-type changing $n$-manifold with a transverse radical, and let $\mathcal{H} \subset \mathbb{R}^{n}$ be a codimension one closed hypersurface of signature change without boundary.\footnote{Here ``closed'' is meant in the topological sense of ``the complement of an open subset of $\mathbb{R}^{n}$'' and not in the manifold sense of ``a manifold without boundary that is compact.''} Then $\mathcal{H}$ is always orientable.
\end{proposition}

\begin{proof}
This can be shown by a purely topological argument, as in~\cite{Samelson }.
\end{proof}

~

\begin{proposition}
Let $(M,g)$ be a transverse, signature-type changing, oriented, $n$-dimensional
manifold with a transverse radical, and let $\mathcal{H}\subset M$
be the hypersurface of signature change. Then $\mathcal{H}$ is also
oriented.
\end{proposition}

\begin{proof}
The hypersurface of signature change, as a closed submanifold of codimension one, is the inverse image of a regular value of a smooth map \( f \colon M \to \mathbb{R} \). Specifically, \(\mathcal{H} = f^{-1}(c)\) for some regular value \(c \in \mathbb{R}\). The manifold \(M\) is oriented, so its tangent bundle \(TM\) is oriented, meaning there is a consistent choice of orientation on each tangent space \(T_p M\) for \(p \in M\). Since \(\mathcal{H}\) is a hypersurface in \(M\), at each point \(q \in \mathcal{H}\), the tangent space \(T_q \mathcal{H}\) is a subspace of the tangent space \(T_q M\) of dimension \(n-1\), and therefore \(T\mathcal{H}\) is a subbundle of \(TM\). The remaining direction in \(T_q M\) can be described by a normal vector \(N(q)\), which is a vector in \(T_q M\) that is perpendicular to \(T_q \mathcal{H}\). 

~

Since \(M\) is oriented, for each point \(q \in \mathcal{H}\), the tangent space \(T_q M\) has an orientation that can be described by an ordered basis, say \(\{v_1, \ldots, v_{n-1}, N(q)\}\), where \(\{v_1, \ldots, v_{n-1}\}\) is an oriented basis for \(T_q \mathcal{H}\) and \(N(q)\) is the normal vector. Hence, this induces a consistent orientation on \(T_q \mathcal{H}\) across all points \(q \in \mathcal{H}\), since the orientation of \(M\) provides a consistent choice of \(N(q)\) across \(\mathcal{H}\). Therefore, \(\mathcal{H}\) inherits a consistent orientation from \(M\), proving that \(\mathcal{H}\) is oriented.

~

Moreover, without loss of generality, we can choose \(1\) as a regular value (see also \cite{Hasse + Rieger}). Thus, \(\mathcal{H} := f^{-1}(1) = \{p \in M \mid f(p) = 1\}\) is a submanifold of \(M\) of dimension \(n-1\). For every \(q \in \mathcal{H}\), the tangent space \(T_q \mathcal{H} = T_q (f^{-1}(1))\) to \(\mathcal{H}\) at \(q\) is the kernel \(\ker(df_q)\) of the map \(df_q \colon T_q M \to T_1 \mathbb{R}\). Then \(T_q \mathcal{H} = \left\langle \text{grad} f_q \right\rangle^\bot\), and therefore the gradient \(\text{grad} f\) yields an orientation of \(\mathcal{H}\).
\end{proof}

~

Provided a transverse, signature-type changing manifold $(M,g)$ with
a transverse radical is pseudo-time orientable, then we can choose
one of the two possible time orientations at any point in each connected
component of $M_{L}$, and thus designating the future direction of
time in the Lorentzian regime. On $\mathcal{H}$ all non-spacelike
vectors are \textit{lightlike} and smoothly divided into two classes
in a natural way: the vectors located at an initial base point on
$\mathcal{H}$ are either pointing towards $M_{L}$ or towards $M_{R}$.
This together with the existent absolute time function (that establishes
a time concept~\cite{Kossowksi + Kriele - Signature type change and absolute time in general relativity}
in the Riemannian region) can be considered as arrow of time on $M$. 

~

\begin{definition}[Natural time
direction]
\label{Definition-(Natural-time direction)} Let $(M,g)$ be a pseudo-time orientable, transverse,
signature-type changing $n$-dimensional manifold with a transverse
radical. Then in the neighborhood of $\mathcal{H}$ the absolute time
function $\mathfrak{h}(t,\mathbf{\hat{x}}):=t$, where $\mathbf{\hat{x}}=(x^{1},\ldots,x^{n-1})$,
imposes \textit{a natural time direction} by postulating that the
future corresponds to the increase of the absolute time function.
In this way, the time orientation is determined in $M_{L}$.
\end{definition}

~

\begin{remark}
\label{rem:future direction}Note that $\partial_{t}$, with an initial
point on $\mathcal{H}$, points in the direction in which $t=\mathfrak{h}(t,\mathbf{\hat{x}})$
increases while $x_{i}$ remains constant. Away from the hypersurface,
the future direction is defined \textit{relative to $\mathcal{H}$}
by the accordant time orientation of $M_{L}$. Recall that functions
of the type, such as the \textit{absolute time function}, typically
lead to metric splittings by default. 
\end{remark}

\begin{figure}[h]
\centering
\includegraphics[width=1.05\textwidth]{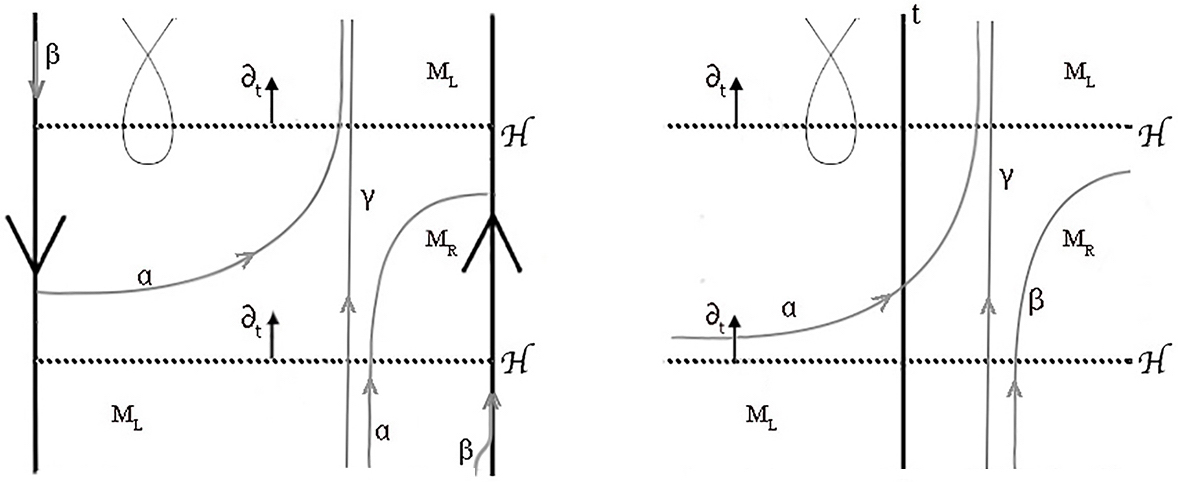}
\caption{In the left example the curves
$\alpha$ and $\gamma$ are both future-directed. The curve $\beta$
runs within the edge that is twisted and identified with the left
edge; therefore $\beta$ is neither future-directed nor past-directed.
In the right example the curves $\alpha$, $\beta$ and $\gamma$
are future-directed. In both examples the loops around $\mathcal{H}$
are neither future-directed nor past-directed.}
\label{fig:Sandwich-structure}
\end{figure}

\begin{definition}[Future-directed]
A pseudo-timelike curve (see Definition~\ref{Definition Pseudo-timelike curve})
in $(M,g)$ is \textit{future-directed} (in the sense of Definition~\ref{Definition-(Natural-time direction)}
and Remark~\ref{rem:future direction}) if for every point in the curve:
\begin{enumerate}[label=(\roman*)]
    \item within $M_{L}$ the tangent vector is future-pointing, and
    \item on $\mathcal{H}$ the associated tangent vector with an initial
    base point on $\mathcal{H}$ is future-pointing.
\end{enumerate}
\end{definition}

\textbf{~}\\ Respective \textit{past-directed curves} are defined
analogously. Notice that, per assumption, one connected component
of $M\setminus\mathcal{H}$ is Riemannian and all other connected
components $(M_{L_{\alpha}})_{\alpha\in I}\subseteq M_{L}\subset M$
are Lorentzian. This configuration could (at least locally) potentially
allow for a $M_{L}-M_{R}-M_{L}$-sandwich-like structure of $M$,
where $\mathcal{H}$ consists of two connected components $(\mathcal{H}_{\alpha})_{\alpha\in\{1,2\}}$.
Consequently, this would also imply the existence of two absolute
time functions, see Figure~\ref{fig:Sandwich-structure}.

~

\begin{definition}[Pseudo-chronological past and future]
Let $(M,g)$ be a
pseudo-time orientable, transverse, signature-type changing $n$-dimensional
manifold with a transverse radical. 

$\mathcal{I}^{-}(p)=\{q\in M\colon q\ll p\}$ is the \textit{pseudo-chronological
past} of the event $p\in M$. In other words, for any two points $q,p\in M$,
we write $q\ll p$ if there is a future-directed pseudo-timelike curve
from $q$ to $p$ in $M$.

$\mathcal{I}^{+}(p)=\{q\in M\colon p\ll q\}$ is the \textit{pseudo-chronological
future} of the event $p\in M$. In other words, for any two points
$p,q\in M$, we write $p\ll q$ if there is a future-directed pseudo-timelike
curve from $p$ to $q$ in $M$.
\end{definition}

\begin{figure}[h]
\centering
\includegraphics[width=0.75\textwidth]{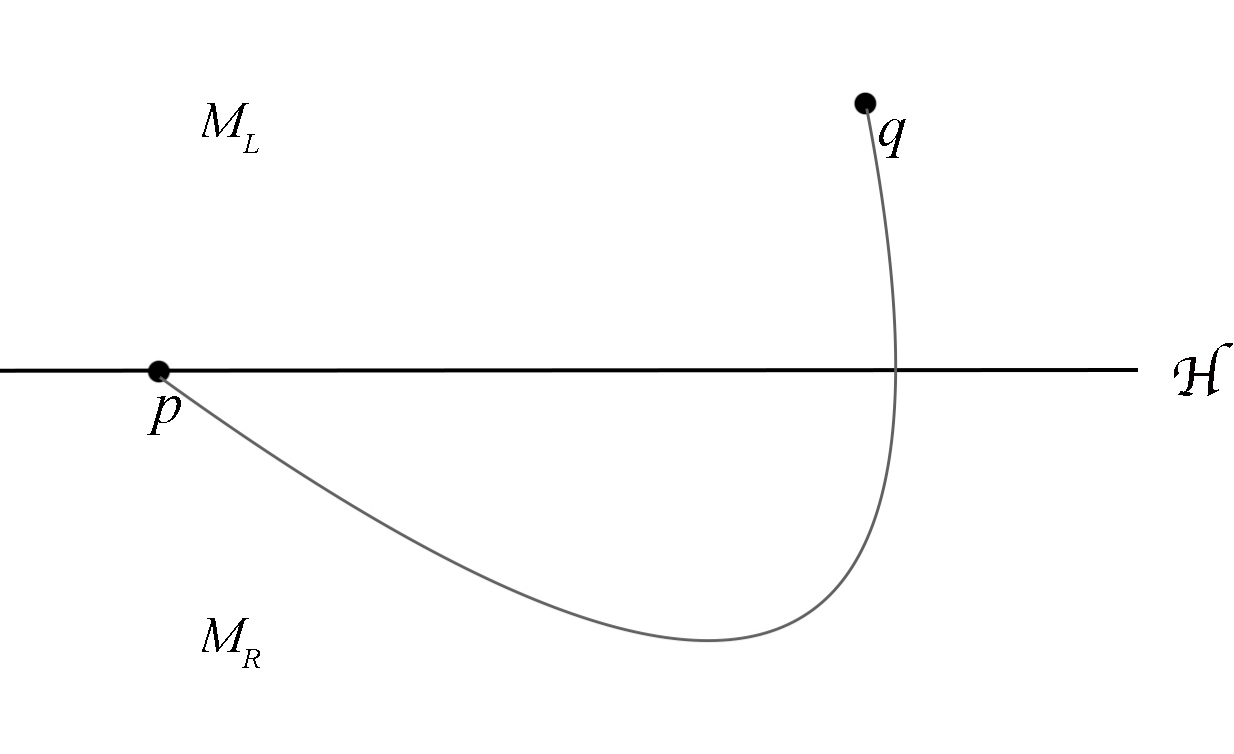}
\caption{For an event $p\in\mathcal{H}$
there exists a future-directed pseudo-timelike curve (as depicted)
that connects the points $p$ and $q$ in $M$. Similarly any point
in $M$ can be reached by such a future-directed pseudo-timelike curve
from $p$. That is why for the pseudo-chronological future we have
$\mathcal{I}^{+}(p)=\{q\in M\colon p\ll q\}=M$.}
\label{fig:chronological future}
\end{figure}

~

\begin{remark}
\label{rem:Pseudo-chronoligical future and past}Interestingly, this
definition leads to the following peculiar situation: Recall that
any curve is denoted pseudo-timelike if its $M_{L}$-segment is timelike.
To that effect, all curves that steer clear of $M_{L}$ (and do not
have a $M_{L}$--segment) are also considered pseudo-timelike. When
$p\in\mathcal{H}\cup M_{R}$ then the pseudo-chronological past of
$p$ is $\mathcal{I}^{-}(p)=M\setminus M_{L}$ and the pseudo-chronological
future of $p$ is $\mathcal{I}^{+}(p)=M$, see Figure~\ref{fig:chronological future}.
\end{remark}

\section{Chronology violating pseudo-timelike loops}

In Section~\ref{sec:Pseudo-Timelike curves}, we introduced the notion
of closed pseudo-timelike curves on a signature-type changing background
and we demonstrated how they must be defined to ensure that the concept
of causality remains meaningful. In this section, we will reveal the
non-well-behaved nature of transverse, signature-type changing, $n$-dimensional
manifolds with a transverse radical. 

\subsection{Local pseudo-timelike loops}\label{Subsec: Local Loops}

In a sufficiently small region near the junction of signature change,
these manifolds exhibit local anomalies. Specifically, each point
on the junction gives rise to the existence of closed time-reversing
loops, challenging conventional notions of temporal consistency. One
of our main results, Theorem~\ref{Theo: Local Loops-1}, can now
be proved quite easily.

~

{
\renewcommand\thetheorem{}
\begin{theorem}[Local loops]
\label{Theo: Local Loops}Let $(M,\tilde{g})$ be a transverse, signature-type
changing, $n$-dimensional ($n\geq2$) manifold with a transverse
radical. Then in each neighborhood of each point ${\color{black}\ensuremath{q\in\mathcal{H}}}$
there always exists a pseudo-timelike loop.
\end{theorem}
}

\begin{proof}
Let $\tilde{g}=-t(dt)^{2}+\tilde{g}_{jk}(t,x^{1},\ldots,x^{n-1})dx^{i}dx^{k}$,
$j,k\in\{1,\ldots,n-1\}$, be a transverse, signature-type changing
metric with respect to a radical-adapted Gauss-like coordinate patch
$(U_{\varphi},\varphi)$ with $U_{\varphi}\cap\mathcal{H}\neq\emptyset$.\footnote{This is, $U_{\varphi}$ is sufficiently small to be expressed in the
adapted radical-adapted Gauss-like coordinate system $\xi(U_{\varphi})$.} Choose smooth coordinates $(t_{0},x_{0}^{1},\ldots,x_{0}^{n-1})$
with $t_{0}>0$ and $\xi_{0}>0$, such that 
\[
C_{0}:=[0,t_{0}]\times B_{\xi_{0}}^{n-1}=[0,t_{0}]\times\{x\in\mathbb{R}^{n-1}\mid{\displaystyle \sum_{k=1}^{n-1}(x^{k})^{2}\leq\xi_{0}^{2}\}\subset}\mathbb{R}^{n}
\]
 is contained in the domain of the coordinate chart (open neighborhood)
$U_{\varphi}$.  Then 
\[
C_{0}\times\mathbb{S}^{n-2}=C_{0}\times\{v\in\mathbb{R}^{n-1}\mid{\displaystyle \sum_{k=1}^{n-1}(v^{k})^{2}=1\}}
\]
as a product of two compact sets is again compact. 

\textbf{~}\\ Next, consider the function

\[
\tilde{G}\colon C_{0}\times\mathbb{S}^{n-2}\longrightarrow\mathbb{R},
\]

\[
(t,x^{1},\ldots,x^{n-1},v^{1},\ldots,v^{n-1})\mapsto\tilde{g}_{jk}(t,x^{1},\ldots,x^{n-1})v^{j}v^{k}.
\]
As $\tilde{G}$ is a smooth function defined on the compact domain
$C_{0}\times\mathbb{S}^{n-2}$, by the Extreme Value Theorem it has
an absolute minimum $G_{0}$. Hence, on $(U_{\varphi},\varphi)$ we
can uniquely define $\tilde{g}_{0}=-t(dt)^{2}+G_{0}\delta_{jk}dx^{j}dx^{k},\:j,k\in\{1,\ldots,n-1\}$.

\textbf{~}\\ By this definition, for all nonzero lightlike vectors
$X\in T_{p}M,\:p\in C_{0}$ with respect to $\tilde{g}_{0}$, we have
$\tilde{g}_{0}=-t(X^{0})^{2}+G_{0}\delta_{jk}X^{j}X^{k}=0\Longleftrightarrow-t(X^{0})^{2}=-G_{0}\delta_{jk}X^{j}X^{k}$,
then

~

\[
\tilde{g}(X,X)=-t(X^{0})^{2}+\tilde{g}_{jk}(t,x^{1},\ldots,x^{n-1})X^{j}X^{k}
\]

\[
=-G_{0}\delta_{jk}X^{j}X^{k}+\tilde{g}_{jk}(t,x^{1},\ldots,x^{n-1})X^{j}X^{k}
\]

\[
=\delta_{jk}X^{j}X^{k}\cdot(-G_{0}+\tilde{g}_{rs}(t,x^{1},\ldots,x^{n-1})\frac{X^{r}}{\sqrt{\delta_{ab}X^{a}X^{b}}}\frac{X^{s}}{\sqrt{\delta_{cd}X^{c}X^{d}}})\geq0.
\]
Clearly, $\tilde{g}(X,X)\geq0$ because $G_{0}>0$ per definition
and $\delta_{jk}X^{j}X^{k}=\frac{t(X^{0})^{2}}{G_{0}}\geq0$. Therefore,
the vector $X\in T_{p}M,\:p\in C_{0}$ is not timelike with respect
to $\tilde{g}$. This means, within $C_{0}$ the $\tilde{g}$-light
cones always reside inside of the $\tilde{g}_{0}$-light cones, i.e.
$\tilde{g}\leq\tilde{g}_{0}$ in $C_{0}$. The cull cones of $\tilde{g}_{0}$
are more opened out than those of the metric $\tilde{g}$. Denote
$p_{0}\in C_{0}$ by $(t(p_{0}),x^{1}(p_{0}),\ldots,x^{n-1}(p_{0}))=(t_{0},x_{0}^{1},\ldots,x_{0}^{n-1})$. 

\textbf{~}\\ As $(M,\tilde{g})$ is an $n$-dimensional manifold
for which in the neighborhood of $\mathcal{H}$ radical-adapted Gauss-like
coordinates exist, we can single out the time coordinate that defines
the smooth absolute time function $t$ whose gradient in $M_{L}$
is everywhere non-zero and timelike. Hence, $(M,\tilde{g})\mid_{U_{\varphi}}$
can be decomposed into spacelike hypersurfaces $\{(U_{\varphi})_{t_{i}}\}$
which are specified as the level sets $(U_{\varphi})_{t_{i}}=t^{-1}(t_{i})$
of the time function.\footnote{This collection of space-like slices $\{(U_{\varphi})_{t}\}$ should
be thought of as a foliation of $U_{\varphi}$ into disjoint $(n-1)$-dimensional
Riemannian manifolds.} The restriction $(\tilde{g}_{0})_{t_{i}}$ of the metric $\tilde{g}_{0}$
to each spacelike slice makes the pair $((U_{\varphi})_{t_{i}},(\tilde{g}_{0})_{t_{i}})$
a Riemannian manifold.

~

\begin{figure}[h]
\centering
\includegraphics[width=0.75\textwidth]{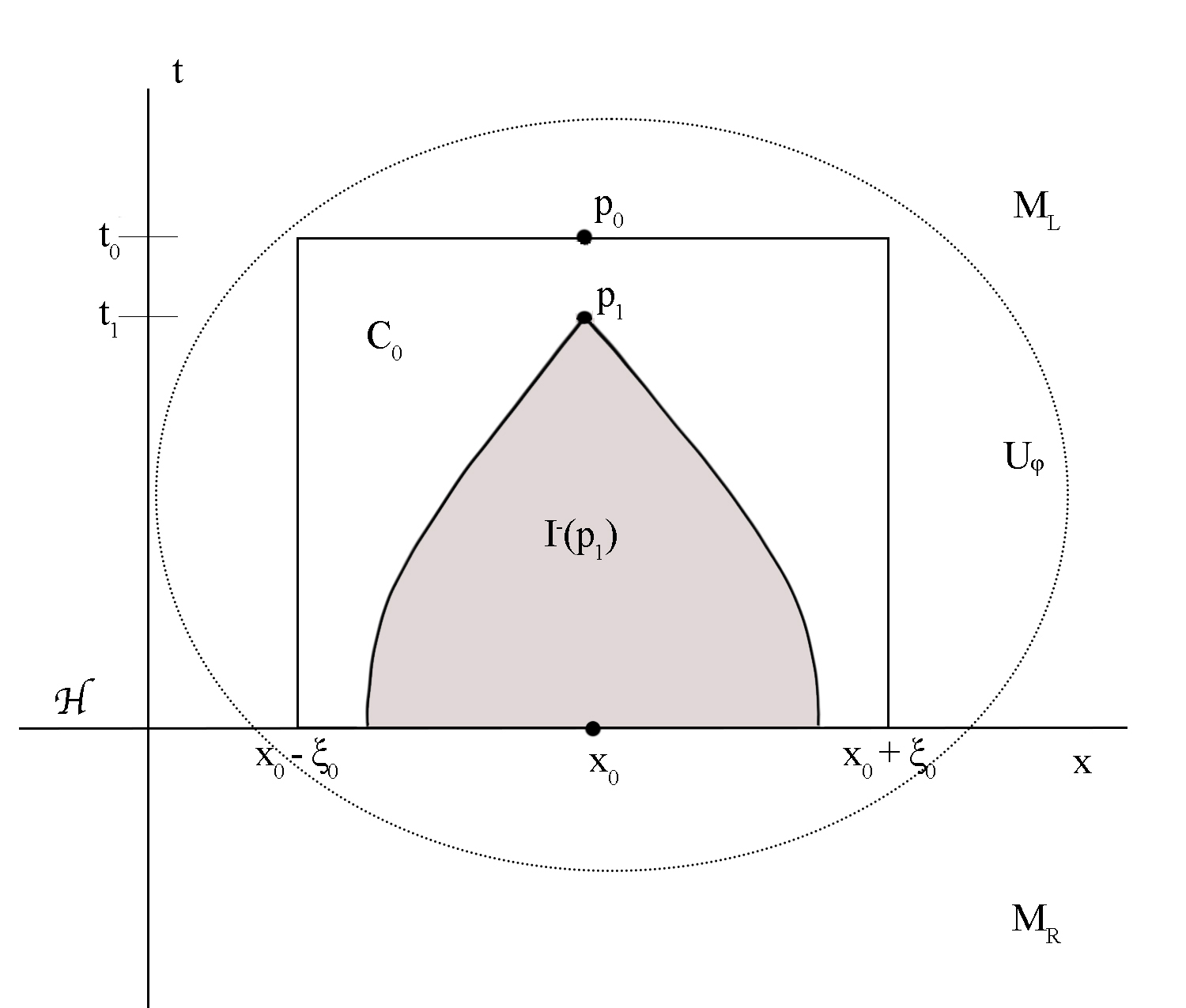}
\caption{The chronological past $I^{-}(p_{1})$
of a point $p_{1}\in U_{\varphi}$.}
\label{fig: Local Loop Theorem}
\end{figure}

For a lightlike curve $\alpha(t)\colon I\longrightarrow U_{\varphi}$
with starting point $p_{0}$, we have $\delta_{jk}\frac{dx^{j}}{dt}\frac{dx^{k}}{dt}>0$
for each slice $(U_{\varphi})_{t_{i}}$ with $t\neq0$. Lightlike
curves with starting point $p_{0}$ can be parametrized with the Euclidean
arc length $\sigma$ in $B_{\xi_{0}}^{n-1}$, such that $(\tilde{g}_{0})_{t}(\dot{\alpha}(\sigma),\dot{\alpha}(\sigma))=\delta_{jk}\frac{dx^{j}}{d\sigma}\frac{dx^{k}}{d\sigma}=1$,
$\forall\:\sigma\in I$, where $I$ is some interval in $\mathbb{R}$.
More precisely, $\sigma$ can be considered as arc length (parameter)
in terms of some auxiliary Riemannian metrics, each defined on a hypersurface
with $t=\textrm{const}$. Consequently we get

\[
0=\tilde{g}_{0}(\dot{\alpha}(\sigma),\dot{\alpha}(\sigma))=-t(\dot{\alpha}^{0})^{2}+G_{0}\delta_{ik}\dot{\alpha}^{j}\dot{\alpha}^{k}
\]

\[
=-t(\frac{dt}{d\sigma})^{2}+G_{0}\underset{1}{\underbrace{\delta_{ik}\frac{dx^{j}}{d\sigma}\frac{dx^{k}}{d\sigma}}=}-t(\frac{dt}{d\sigma})^{2}+G_{0},
\]
and this implies 
\[
\frac{d\sigma}{dt}=\pm\sqrt{\frac{t}{G_{0}}}\Longrightarrow\sigma(t)=\pm\int\sqrt{\frac{t}{G_{0}}}dt=\pm\frac{2}{3}t\sqrt{\frac{t}{G_{0}}}+\textrm{const}.
\]

\textbf{~}\\ Since $\sigma$ is given as a function of $t$, it represents
the arc length from the starting point at $t(p_{0})=t_{0}$ to $t(0)=0$.
Then past-directed $\tilde{g}_{0}$-lightlike curves emanating from
$p_{0}$ reach the hypersurface at $t=0$ after passing through the
arc length distance 
\[
\triangle\sigma=\pm{\displaystyle \int_{0}^{t_{0}}}\sqrt{\frac{t}{G_{0}}}dt=\pm\frac{2}{3}t_{0}\sqrt{\frac{t_{0}}{G_{0}}}+\textrm{const}.=\pm\frac{2}{3}\sqrt{\frac{t_{0}^{3}}{G_{0}}}+\textrm{const}
\]
along the said section of the curve from the fixed starting point
$p_{0}$.

\textbf{~}\\ Provided this arc length distance satisfies $\triangle\sigma\leq\xi_{0}$,
then the past-directed lightlike curves $\alpha(t)$ (emanating from
$p_{0}$) reach the hypersurface at $t=0$ while remaining within
$C_{0}$. Accordingly this is also the case for $\tilde{g}$-lightlike
curves emanating from $p_{0}$. Conversely, if $\triangle\sigma>\xi_{0}$
then there exist past-directed $\tilde{g}_{0}$-lightlike curves emanating
from $p_{0}$ that reach the hypersurface outside of $C_{0}$.\textbf{~}\\ In
this case, we have $\triangle\sigma=\frac{2}{3}\sqrt{\frac{t_{0}^{3}}{G_{0}}}>\xi_{0}\Longleftrightarrow t_{0}>\sqrt[3]{\frac{9}{4}\xi_{0}^{2}\cdot G_{0}}$
and we must adjust the new starting point $p_{1}=(t_{1},x_{0})$ accordingly
by setting $t_{1}\leq\sqrt[3]{\frac{9}{4}\xi_{0}^{2}\cdot G_{0}}<t_{0}$.
Thereby we make sure that all past-directed $\tilde{g}$-lightlike
curves emanating from $p_{1}$ hit the hypersurface $\mathcal{H}$
without leaving $C_{0}$. That is $\overline{I_{\tilde{g}}^{-}(p_{1})}\subset C_{0}\subset U_{\varphi}\subset M$,
where $\overline{I_{\tilde{g}}^{-}(p_{1})}$ is the \textit{$\tilde{g}$-chronological
past} of the event $p_{1}\in M_{L}$, restricted to $M_{L}\cup\mathcal{H}$,
see Figure~\ref{fig: Local Loop Theorem}.

\textbf{~}\\ It now suffices to connect two of such points $\hat{x}_{1},\hat{x}_{2}\in\overline{I_{\tilde{g}}^{-}(p_{0})}\cap\mathcal{H}$
(or, if need be $\overline{I_{\tilde{g}}^{-}(p_{1})}\cap\mathcal{H}$)
in an arbitrary fashion within the Riemannian sector $M_{R}$. By
what a pseudo-timelike loop gets generated, if $U_{\varphi}$ was
chosen small enough.
\end{proof}

~

Summarized, for each neighborhood $U(q)$ that admits radical-adapted
Gauss-like coordinates $\xi=(t,\mathbf{\hat{x}})=(t,x^{1},\ldots,x^{n-1})$
centered at some $q\in\mathcal{H}$, and $U(q)\cap\mathcal{H}\neq\emptyset$,
we are able to pick a point $p_{0}\in U(q)$ and an associated compact
set $C_{0}\subset U(q)$. For the metric $\tilde{g}$ there exists
a corresponding uniquely (i.e., only dependent on the chosen set $C_{0}$)
defined metric $\tilde{g}_{0}$ with $\tilde{g}\leq\tilde{g}_{0}$
within $C_{0}$.\footnote{The set \( C_{0} \) does not need to be ``maximal'' (in some sense) and is therefore not unique.} Then we must distinguish between two cases, that is

~

i)$\:$ with respect to the metric $\tilde{g}_{0}$ we have $I_{0}^{-}(p_{0})\subset C_{0}$,
then also $I^{-}(p_{0})\subset C_{0}$ with respect to $\tilde{g}$,

ii) with respect to the metric $\tilde{g}_{0}$ we have the situation
$I_{0}^{-}(p_{0})\nsubseteq C_{0}$, then there exists a point $p_{1}=(t_{1},x_{0})\in C_{0}\setminus\mathcal{H}$
with $t_{1}<t_{0}$, such that $I_{0}^{-}(p_{1})\subset C_{0}$, hence
also $I^{-}(p_{1})\subset C_{0}$ with respect to $\tilde{g}$.

~

Thus, for any point $q\in\mathcal{H}$ we can find a sufficiently
small neighborhood $\tilde{U}\subset U(q)$ containing a point $p\in M_{L}$,
such that all past-directed, causal curves emanating from that point,
reach the hypersurface within a sufficiently small set $C_{0}$.

~

\begin{corollary}
Let $(M,\tilde{g})$ be a transverse, signature-type changing, $n$--dimensional
manifold with a transverse radical. Then in each neighborhood of each
point ${\color{black}\ensuremath{q\in\mathcal{H}}}$ there always
exists a pseudo-lightlike curve.
\end{corollary}

~

The above corollary follows directly from Theorem~\ref{Theo: Local Loops-1}
because we have proven that all past-directed $\tilde{g}$-\textit{lightlike}
curves emanating from $p_{1}$ hit the hypersurface $\mathcal{H}$
without leaving $C_{0}$. That is $\overline{I_{\tilde{g}}^{-}(p_{1})}\subset C_{0}\subset U_{\varphi}\subset M$,
where $\overline{I_{\tilde{g}}^{-}(p_{1})}$ is the closure of the
\textit{$\tilde{g}$-chronological past} of the event $p_{1}\in M_{L}$.
Hence, we have also shown that the causal past is within
$C_{0}$, and furthermore, $J_{\tilde{g}}^{-}(p_{0})\subset\overline{I_{\tilde{g}}^{-}(p_{0})}$. 

And since in every neighborhood of each point $q\in\mathcal{H}$ there
always exists a pseudo-timelike loop, we can straightforwardly assert
the following

~

\begin{corollary}
A transverse, signature-type changing manifold $(M,\tilde{g})$ with
a transverse radical has always time-reversing pseudo-timelike loops.
\end{corollary}

~

\begin{remark}We note that the loops constructed in the proof of Theorem~\ref{Theo: Local Loops-1} are entirely contained within the small neighborhoods $U(q)$ and do not extend globally. More generally, pseudo-timelike loops may exist that traverse larger portions of the manifold, particularly if the hypersurface $\mathcal{H}$ has a toroidal or otherwise nontrivial global structure. Allowing self-intersections does not conflict with treating the objects as parameterized curves.
\end{remark}
~

As a matter of course, in the Lorentizan region the tangent space
at each point is isometric to Minkowski space which is time orientable.
Hence, a Lorentzian manifold is always infinitesimally time- and space-orientable,
and a continuous designation of future-directed and past-directed
for non-spacelike vectors can be made (infinitesimally and therefore,
by continuity, also locally).\footnote{In case the Lorentzian manifold is time-orientable, a continuous designation
of future-directed and past-directed for non-spacelike vectors can
be made allover.}

\textbf{~}\\ Having said that, the infinitesimal properties of a
manifold with a signature change are identical to those of a Lorentzian
manifold only within the Lorentzian sector. However, when examining
the Riemannian sector and the hypersurface, specific distinctions
arise. The Riemannian sector and the hypersurface are not infinitesimally
modelable by a Minkowski space. While the Riemannian sector reveals
an absence of a meaningful differentiation between past- and future-directed
vectors, on the hypersurface, one has the flexibility to make arbitrary
assignments of such distinctions at the infinitesimal level. If one
now determines on the hypersurface whether the direction towards the
Lorentzian sector is the future or past direction, it is not only
a reference to the tangent space at a point. Rather, it is a local
consideration.

\textbf{~}\\ In the context of local considerations, in a Lorentzian
manifold the existence of a timelike loop that flips its time orientation
(i.e. the timelike tangent vector switches between the two designated
components of the light cone) is a sufficient condition for the absence
of time orientability. Based on the previous theorem (at the beginning
of the present subsection), this is also true for a transverse, signature-type
changing manifold $(M,\tilde{g})$ with a transverse radical: As we
have proved above, through each point on the hypersurface $\mathcal{H}$
we have locally a closed time-reversing loop. That is, there always
exists a closed pseudo-timelike path in $M$ around which the direction
of time reverses, and along which a consistent designation of future-directed
and past-directed vectors cannot be defined. 

~

\begin{figure}[h]
\centering
\includegraphics[width=0.55\textwidth]{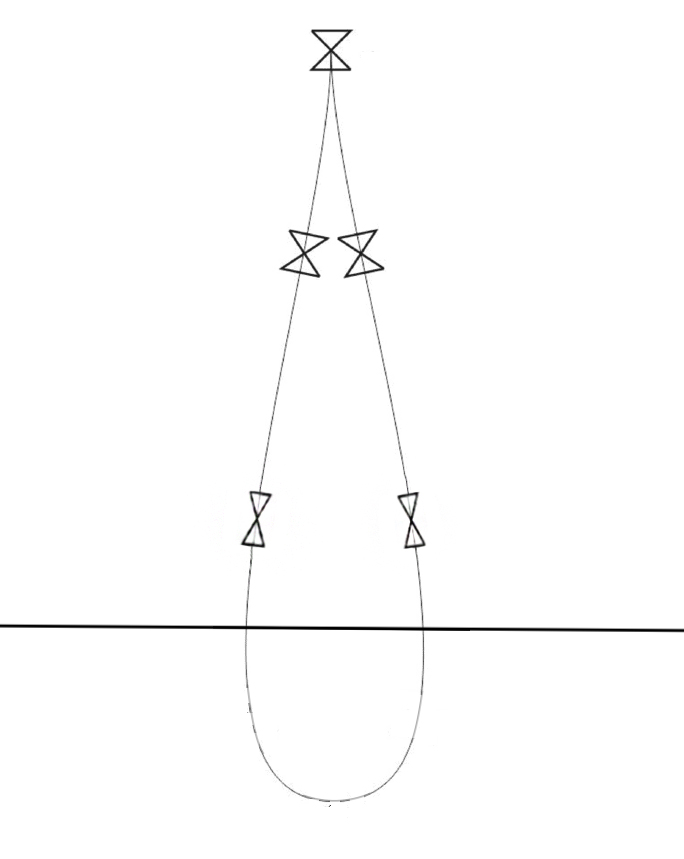}
\caption{A closed time-reversing loop.}
\label{fig:A-closed-time-reversing-loop}
\end{figure}

An observer in the region $M_{L}$ near $\mathcal{H}$ perceives these
locally closed time-reversing loops (Figure~\ref{fig:A-closed-time-reversing-loop})
as the creation of a particle and an antiparticle at two different
points $\hat{q},q\in\mathcal{H}$.\footnote{Such locally closed time-reversing loops around $\mathcal{H}$ obviously
do not satisfy the causal relation $\ll$ as introduced above.} This could be taken as an object entering the Riemannian region,
then resurfacing in the Lorentzian region and proceeding to propagate “forward in time” with respect to the time orientation in $M_L$. Since the mathematical definition of the loop itself does not privilege a direction of traversal, this picture is symmetric.

~

So in a transverse, signature-type changing manifold $(M,\tilde{g})$,
the hypersurface with its time-reversing loops could be tantamount
to a region of particle-antiparticle origination incidents. Moreover,
Hadley~\cite{Hadley - The orientability of spacetime} shows for
Lorentzian spacetimes that a failure of time-orientability of a spacetime
region is indistinguishable from a particle-antiparticle annihilation
event. These are then considered equivalent descriptions of the same
phenomena. It would be interesting to explore how this interpretation
can be carried over to signature-type changing manifolds. \textbf{~}\\{}

For fields, take the conjugate \(\psi_{t}^{A} = e^{-i\hat{H}t}\psi^{*}\) of \(\psi_{t} = e^{i\hat{H}t}\psi\): The unitary temporal evolution of the field operator for antiparticles arises from the temporal evolution of the field operator for particles by applying the same Hamiltonian operator to the adjoint field operator under time reversal. Some literature~\cite{Hadley - The Logic of Quantum Mechanics Derived from Classical General Relativity} points to the idea that concepts in quantum field theory are predicated on acausal properties derived from general relativity. In this context, Blum et al.~\cite{Blum et al} stress the importance of the CPT theorem (quoting verbatim):

\begin{quote}
``CPT theorem is the statement that nothing would change---nobody would notice and the predictions of physics would not be altered---if we simultaneously replace particles by antiparticles and vice versa. Replace everything by its mirror image or, more exactly, exchange left and right, up and down, and front and back, and reverse the flow of time. We call this simultaneous transformation CPT, where C stands for Charge Conjugation (exchanging particles and antiparticles), P stands for parity (mirroring), and T stands for time reversal.''
\end{quote}

\subsection{Global pseudo-timelike loops}

The existence of such pseudo-timelike curves locally near the hypersurface
that loop back to themselves, gives rise to the question whether this
type of curves also occur globally. We want to elucidate this question
in the following.\footnote{A spacetime is a Lorentzian manifold that models space and time in
general relativity and physics. This is conventionally formalized
by saying that a spacetime is a smooth connected time-orientable Lorentzian
manifold $(M,g)$ with $\dim M=4$. But in what follows we want to
study the $n$-dimensional ($n\geq2$) case.}

~

\begin{definition}[Stably causal]
~\cite{Minguzzi} A connected time-orientable
Lorentzian manifold $(M,g)$ is said to be \textit{stably causal}
if there exists a nowhere-vanishing timelike vector field $V_{a}$
such that the Lorentzian metric on $M$ given by $g':=g_{ab}-V_{a}V_{b}$
admits no closed causal curves. In other words, if $(M,g)$ is stably
causal then, for some timelike $V_{a}$, the metric $g':=g_{ab}-V_{a}V_{b}$
on $M$ is causal.
\end{definition}

~

\begin{remark}
A partial ordering $<$ is defined in the set of all Lorentzian metrics
$Lor(M)$ on $M$ in the following way: $g<g'$ iff all causal vectors
for $g$ are timelike for $g'$. Then the metric $g_{\lambda}=g+\lambda(g'-g)$,
$\forall\:\lambda\in[0,1]$ is a Lorentzian metric on $M$, as well.
Also, recall that $g<g'$ means that the causal cones of $g$ are
contained in the timelike cones of $g'$. A connected time-orientable
Lorentzian manifold $(M,g)$ is stably causal if there exists $g'\in Lor(M)$,
such that $g'>g$, with $g'$ causal.
\end{remark}

~

\begin{lemma}
~\textup{\cite{Sanchez} }Stable causality is the necessary and sufficient
condition for the existence of a smooth global time function, i.e.
a differentiable map $T\colon M\rightarrow\mathbb{R}$ such that whenever
$p<<q$ $\Longrightarrow T(p)<T(q)$.
\end{lemma}

~

\begin{definition}[Globally hyperbolic]
~\cite{Bernal + Sanchez - Further Results on the Smoothability of Cauchy Hypersurfaces and Cauchy Time Functions,Hawking + Ellis - The large scale structure of spacetime}
A connected, time-orientable Lorentzian manifold $(M,g)$ is called
\textit{globally hyperbolic} if and only if it is diamond-compact
and causal, i.e., $p\notin J^{+}(p)\:\forall p\in M$.\footnote{Diamond-compact means $J(p,q):= J^{+}(p)\cap J^{-}(q)$ is
compact for all $p,q\in M$. Note that $J(p,q)$ is possibly empty.}
\end{definition}

~

An equivalent condition for global hyperbolicity is as follows~\cite{Geroch - Domain of Dependence}.\\
\begin{definition}
A connected, time-orientable Lorentzian manifold $(M,g)$ is called
\textit{globally hyperbolic} if and only if $M$ contains a Cauchy
surface. A Cauchy hypersurface in $M$ is a subset $S$ that is intersected
exactly once by every inextensible timelike curve in $M$.\footnote{An inextensible curve is a general term that refers to a curve with
no endpoints; it either extends infinitely or it closes in on itself
to form a circle---a closed curve. Specifically, an inextensible
timelike curve is a curve that remains timelike throughout its entire
length and cannot be extended further within the spacetime. In mathematical
terms, a map $\alpha\colon(a,b)\rightarrow M$ is an inextensible
timelike curve in $(M,g)$ if $\alpha(t)$ does not approach a limit
as $t$ increases to $b$ or decreases to $a$, and $\alpha(t)$ remains
timelike for all $t\in(a,b)$. This distinguishes it from inextensible
curves of other causal types, such as null or spacelike curves.}
\end{definition}

~

In 2003, Bernal and S\'{a}nchez~\cite{Bernal+Sanchez - On smooth Cauchy hypersurfaces and Gerochs splitting theorem}
showed that any globally hyperbolic Lorentzian manifold $M$ admits
a smooth spacelike Cauchy hypersurface $S$, and thus is diffeomorphic
to the product of this Cauchy surface with $\mathbb{R}$, i.e. $M$
splits topologically as the product $\mathbb{R}\times S$. Specifically,
a globally hyperbolic manifold is foliated by Cauchy surfaces.\\

\begin{remark}
If $M$ is a smooth, connected time-orientable Lorentzian manifold
with boundary, then we say it is globally hyperbolic if its interior
is globally hyperbolic. 
\end{remark}

~

The next theorem (Global Loops Theorem) is partially based on the
Local Loops Theorem~\ref{Theo: Local Loops-1} and can be considered
a generalization to the global case.\\
{
\renewcommand\thetheorem{}
\begin{theorem}[Global loops]
\label{thm:Global Loop-1} Let $(M,\tilde{g})$ be a pseudo-time orientable,
transverse, signature-type changing, $n$-dimensional ($n\geq2$)
manifold with a transverse radical, where $M_{L}=M\setminus(M_{R}\cup\mathcal{H})$
is globally hyperbolic. Assume that a Cauchy surface $S$ is a subset
of the neighborhood $U=\bigcup_{q\in\mathcal{H}}U(q)$ of $\mathcal{H}$,
i.e. $S\subseteq(U\cap M_{L})=\bigcup_{q\in\mathcal{H}}(U(q)\cap M_{L})$,
with $U(q)$ being constructed as in Theorem~\ref{Theo: Local Loops-1}.
Then for every point $p\in M$, there exists a pseudo-timelike loop
such that $p$ is a point of self-intersection.
\end{theorem}
}

\begin{proof}
Let $(M,\tilde{g})$ be a \textit{pseudo-time orientable} transverse,
signature-type changing, $n$-dimensional ($n\geq2$) manifold with
a transverse radical, where $M_{L}$ is globally hyperbolic with $\tilde{g}\mid_{M_{L}}=g$.
Moreover, there is a neighborhood $U=\bigcup_{q\in\mathcal{H}}U(q)$
of $\mathcal{H}$ sufficiently small to satisfy the conditions for
Theorem~\ref{Theo: Local Loops-1}, and per assumption there exists
a Cauchy surface $S_{\varepsilon}\subseteq(U\cap M_{L})$, $\varepsilon>0$.

\textbf{~}\\ Due to~\cite{Bernal+Sanchez - On smooth Cauchy hypersurfaces and Gerochs splitting theorem}
we know that $M_{L}$ admits a splitting $M_{L}=(\mathbb{R}_{>0})_{t}\times S_{t}=\bigcup_{t\in\mathcal{\mathbb{R}}_{>0}}S_{t}$,
such that the Lorentzian sector $M_{L}$ is decomposed into hypersurfaces
(of dimension $n-1$), specified as the level surfaces $S_{t}=\mathcal{T}^{-1}(t)=\{p\in M_{L}\colon\mathcal{T}(p)=t\},t\in\mathbb{R}_{>0}$,
of the real-valued smooth temporal function $\mathcal{T}\colon M_{L}\longrightarrow\mathbb{R}_{>0}$
whose gradient $\textrm{grad}\mathcal{T}$ is everywhere non-zero
and, clearly, $d\mathcal{T}$ is an exact $1$-form. Within the neighborhood
$U=\bigcup_{q\in\mathcal{H}}U(q)$ this foliation $\bigcup_{t\in\mathcal{\mathbb{R}}_{>0}}S_{t}$
can be chosen in such a way that it agrees with the natural foliation
given by the absolute time function $\mathfrak{h}(t,\mathbf{\hat{x}}):=t$, see
Remark~\ref{rem:future direction} and Definition~\ref{Definition-(Natural-time direction)}.\footnote{Recall that a smooth function $T\colon M\longrightarrow\mathbb{R}$
on a connected time-orientable Lorentzian manifold $(M,g)$ is a global
time function if $T$ is strictly increasing along each future-pointing
non-spacelike curve. Moreover, a temporal function is a time function
$T$ with a timelike gradient $\textrm{grad}T$ everywhere.

Since $M_{L}$ is globally hyperbolic it admits a smooth global time
function $T$ and consequently it admits~\cite{Minguzzi} a temporal
function $\mathcal{T}$. Hence, in the Lorentzian sector $M_{L}$
there exists a global temporal function $\mathcal{T}\colon M_{L}\longrightarrow\mathbb{R}_{>0}$,
and $\textrm{grad}\mathcal{T}$ is orthogonal to each of the level
surfaces $S_{t}=\mathcal{T}^{-1}(t)=\{p\in M_{L}\colon\mathcal{T}(p)=t\},t\in\mathbb{R}_{>0}$,
of $\mathcal{T}$. Note that $\mathcal{T}=t$ is a scalar field on
$M_{L}$, hence $\textrm{grad}\mathcal{T}=\textrm{grad}t=(dt)^{\#}$.}

\textbf{~}\\ Moreover, the level surfaces $(S_{t})_{t\in\mathbb{R}}$
are Cauchy surfaces and, accordingly, each inextensible pseudo-timelike
curve in $M_{L}$ can intersect each level set $S_{t}$ exactly once
as $\mathcal{T}$ is strictly increasing along any future-pointing
pseudo-timelike curve.\footnote{Since $\mathcal{T}$ is regular the hypersurfaces $S_{t}$ never intersect,
i.e. $S_{t}\cap S_{t'}=\emptyset$ for $t\neq t'$.} Then, these level-sets $S_{t}$ are all space-like hypersurfaces
which are orthogonal to a timelike and future-directed unit normal
vector field $n$.\footnote{In other words, the unit vector $n$ is normal to each slice $S_{t}$,
and $g$ restricted to $S_{t}$ is Riemannian.}

\textbf{~}\\ For $\varepsilon$ sufficiently small, the level Cauchy
surface 
\[
S_{\varepsilon}=\mathcal{T}^{-1}(\varepsilon)=\{p\in M_{L}\colon\mathcal{T}(p)=\varepsilon\},\varepsilon\in\mathbb{R}_{>0}
\]
 is contained in $U\cap M_{L}=\bigcup_{q\in\mathcal{H}}(U(q)\cap M_{L})$.\footnote{This is true because all neighborhoods \( U(q) \) with \( q \in \mathcal{H} \)
can be chosen such that the sets \( U(q) \) have a compact closure.
Thus, the \( \overline{U(q)} \) are not ``infinitely wide'', and there exists a strictly positive value \( \varepsilon_{\max} \),
such that for all \( \varepsilon < \varepsilon_{\max} \), the level Cauchy
surface \( S_{\varepsilon} \) is contained in \( U \cap M_{L} \).}

\textbf{~}\\ Therefore, based on Theorem~\ref{Theo: Local Loops-1},
for any $p=(\varepsilon,\mathbf{\hat{x}})\in S_{\varepsilon}\subseteq$$(U\cap M_{L})$
all past-directed and causal curves emanating from that point reach
the hypersurface $\mathcal{H}$. The global hyperbolicity of $M_{L}$
implies that every non-spacelike curve in $M_{L}$ meets each $S_{t}$
once and exactly once since $S_{t}$ is a Cauchy surface. In particular,
the spacelike hypersurface $S_{\varepsilon}$ is a Cauchy surface
in the sense that for any $\bar{p}\in M_{L}$ in the future of $S_{\varepsilon}$,
all past pseudo-timelike curves from $\bar{p}$ intersect $S_{\varepsilon}$.
The same holds for all future directed pseudo-timelike curves from
any point in $M_{L}$ in the past of $S_{\varepsilon}$.

~

\begin{figure}[h]
\centering
\includegraphics[width=0.7\textwidth]{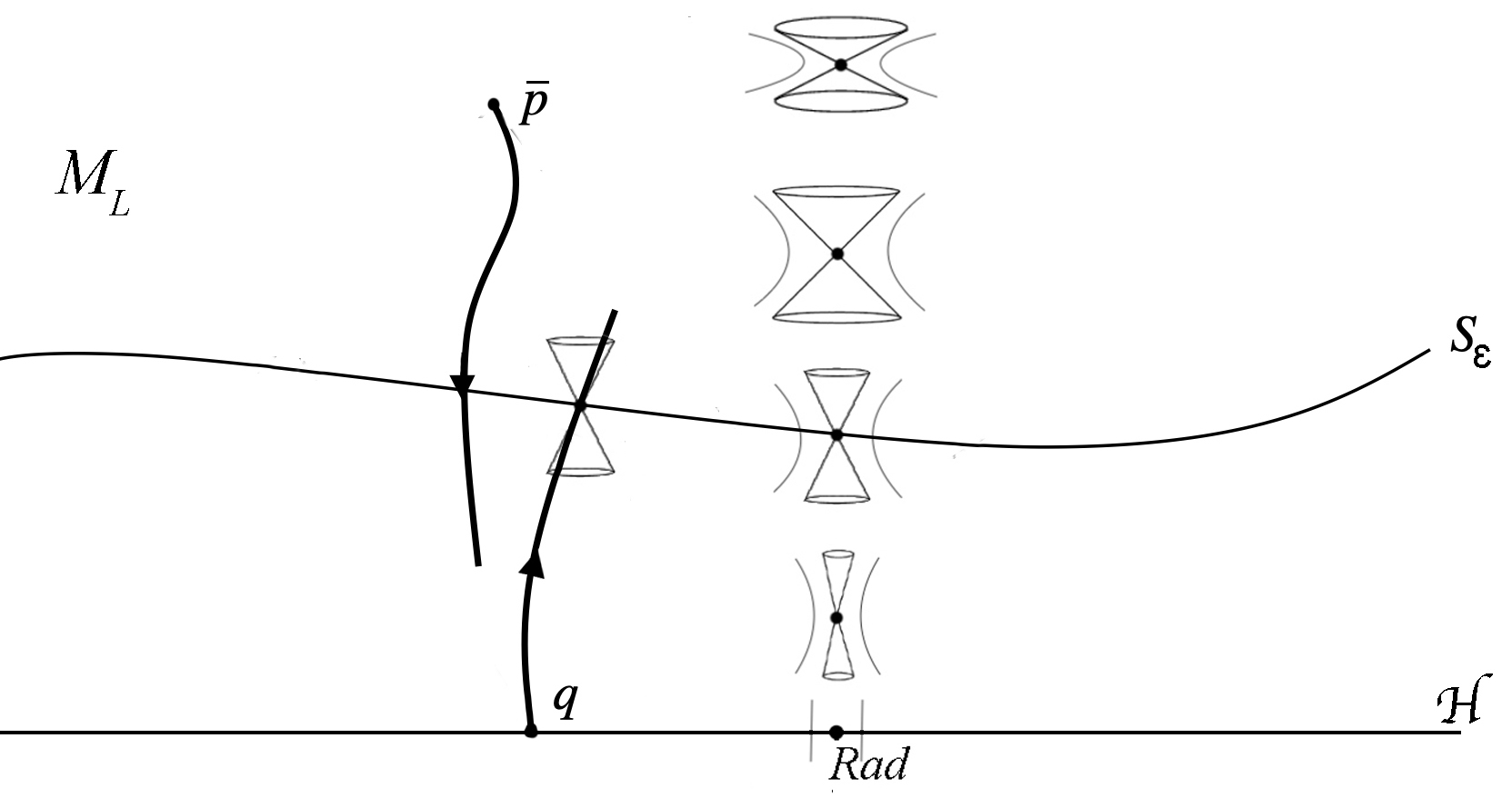}
\caption{In the future of the Cauchy surface $S_{\varepsilon}$, any point 
$\bar{p} \in M_{L}$ with $t > \varepsilon$ has all of its past-directed 
pseudo-timelike curves intersect $S_{\varepsilon}$. Similarly, all 
future-directed pseudo-timelike curves from a point $q \in \mathcal{H}$ 
intersect $S_{\varepsilon}$. This means there must be a point 
$\bar{q} \in \mathcal{H}$ from which $\bar{p} \in M_{L}$ can be reached 
by a future-directed pseudo-timelike curve.}
\label{fig:Cauchy-surface}
\end{figure}

Consequently, by virtue of Theorem~\ref{Theo: Local Loops-1} and
the above argument, all past-directed pseudo-timelike curves emanating
from any $\bar{p}\in M_{L}$ \textit{reach} the hypersurface $\mathcal{H}$.
Analogously we can conclude that any point $\bar{p}\in M_{L}$ \textit{can
be reached} by a future-directed pseudo-timelike curve starting at
some suitable point in $\mathcal{H}$. Recall that, based on Remark~\ref{rem:Pseudo-chronoligical future and past},
we also know that $\mathcal{I}^{+}(q)=\{p\in M:q\ll p\}=M$, that
is, any point in ${\color{black}M=M_{R}\cup\mathcal{H}\cup M_{L}}$
can be reached by a future-directed pseudo-timelike curve from $q\in\mathcal{H}$,
see Figure~\ref{fig:Cauchy-surface}.

\textbf{~}\\ We now obtain a loop with intersection point $p$ in
$M_{L}$ if, for sufficiently small $\varepsilon$, we first prescribe
the intersection point $p=(\varepsilon,\mathbf{\hat{x}})\in S_{\varepsilon}$.
And then we connect the two points lying in $\mathcal{H}$ of the
intersecting curve sections through an arbitrary curve segment in
the Riemannian sector $M_{R}$ (through a suitable choice of the two
curve segments, we can ensure that different points on $\mathcal{H}$
are obtained).

\textbf{~}\\This reasoning also applies to points located on the hypersurface or within 
the Riemannian region. In these cases, the situation is as follows:\\
\begin{enumerate}
    \item If the given point lies on the hypersurface, $p \in \mathcal{H}$, choose a 
    timelike curve segment in $M_L$ connecting $p$ to $S_\varepsilon$ (with 
    $\varepsilon$ sufficiently small). From there, continue along a second timelike curve segment back to another point on $\mathcal{H}$. Finally, connect both hypersurface points through the Riemannian sector. At $p$, the Riemannian curve segment must form a kink where it meets the timelike segment, and the two intersecting curve segments must be extended so that together they produce a self-intersection exactly at $p$.
    
    \item If the given point lies in the Riemannian sector, \(p \in M_R\), begin with a loop of the type constructed in the proof of Theorem~\ref{Theo: Local Loops-1}. Modify the Riemannian portion of this loop so that it passes through \(p\) and forms a self-intersection there. As in the previous case, one intersection point of the loop necessarily lies in \(M_L\), while the second intersection is realized at \(p\).  
\end{enumerate}
\end{proof}

\begin{remark}
Theorem~\ref{Theo: Global Loops-1} explicitly states that through every point in $M$, there always exists a pseudo-timelike loop. Consequently, this applies equally to points on the hypersurface as well as to those in the Riemannian region. In such cases, the construction necessarily produces a pseudo-timelike loop with two intersection points. Indeed, the point $p$ is required to serve as a point of self-intersection. Consequently, the pseudo-timelike loop must also admit an intersection in $M_L$, and at least one of the intersection points must lie in $M_L$ in order for the curve to qualify as a pseudo-timelike loop. The second intersection point is then determined, either as $p$ on $\mathcal{H}$ or as $p$ in $M_R$.\footnote{The construction for points $p \in \mathcal{H}$ necessitates that the timelike curve segments from $M_L$ approach the ``pinched''  light cones near the hypersurface. Consequently, the curve segments that pass through $\mathcal{H}$ into the Riemannian sector $M_R$ and self-intersect at $p \in \mathcal{H}$ must meet with first-order contact.} Topologically, a loop with two intersection points falls within the class of curves described in Definition~\ref{Definition Chronology violating curve}.
\end{remark}

~

\begin{example}
\textit{\label{Example-17}} The prototype of a spacetime $M$ with
signature-type change is obtained by cutting an $S^{4}$ along its
equator and joining it to the corresponding half of a de Sitter space.
It is a well-known fact that the full de Sitter spacetime is globally
hyperbolic~\cite{Kroon. Anti-de Sitter-like spacetimes}, with the
entire manifold possessing a Cauchy surface. When we restrict to half
de Sitter space---by choosing an appropriate region bounded by a
Cauchy surface---this region retains global hyperbolicity. This is
because the Cauchy surface of the full de Sitter spacetime remains
valid in the half-space, ensuring that every inextensible non-spacelike
curve still intersects this surface exactly once. As a result, the
Lorentzian sector, which corresponds to half de Sitter space, is also
globally hyperbolic. Consequently, there are chronology-violating
pseudo-timelike loops through each point in $M$.
\end{example}

~

\begin{corollary}
Let $(M,\tilde{g})$ be a pseudo-time orientable, transverse, signature-type
changing, $n$-dimensional ($n\geq2$) manifold with a transverse
radical, where $M_{L}$ is globally hyperbolic, and ${\color{black}S\subseteq(U\cap M_{L})=\bigcup_{q\in\mathcal{H}}(U(q)\cap M_{L})}$
for a Cauchy surface $S$. Then through every point there exists a
path on which a pseudo-time orientation cannot be defined.
\end{corollary}

\section{Final thoughts}

The intriguing facet of the potential existence of closed timelike
curves within the framework of Einstein's theory lies in the physical
interpretation that CTCs, serving as the worldlines of observers,
fundamentally permit an influence on the causal past. This can also
be facilitated through a causal curve in the form of a loop, i.e.,
the curve intersects itself. In the case of a non-time-orientable
manifold, there would then be the possibility that at the intersection,
the two tangent vectors lie in different components of the light cone.
Thus, the ``time traveler'' at the encounter
with himself, which he experiences twice, may notice a reversal of
the past and future time directions in his surroundings during the
second occurrence, even including the behavior of his younger
version. Regardless of whether this effect exists or not, during the
second experience of the encounter, which he perceives as an encounter
with a younger version of himself, the traveler can causally influence
this younger version and its surroundings.

~

It is important to emphasize that the existence of pseudo-timelike loops near the hypersurface of signature change, as guaranteed by the Loop Theorems, does not imply a breakdown of physical causality. Locally in the Lorentzian region $M_L$, the directions of certain classes of curves are constrained, which can be interpreted physically as limitations on the propagation of signals or particle interactions along these curves. In contrast, in the Riemannian region $M_R$, two competing interpretations are possible regarding the physical meaning of curve segments:

~

\begin{enumerate}[label=(\arabic*)]
    \item \textbf{Unrestricted propagation interpretation:} While $M_R$ allows curves mathematically without directional constraints, one can treat them as a continuation of possible influence from $M_L$. Local causality in $M_L$ remains intact: in a sufficiently small neighborhood of a point $p$, no point $q$ can both influence and be influenced by $p$. In this reading, pseudo-timelike loops are mathematically extended curves that do not correspond to actual causal violation.\\
    
    \item \textbf{Causal barrier interpretation:} Alternatively, $M_R$ can be viewed as a region where no physical propagation occurs. Curves there are spacelike in the Lorentzian sense and cannot carry influence. Consequently, any pseudo-timelike loop passing through $M_R$ cannot induce causality violations even in a global sense.
\end{enumerate}

\textbf{~}\\
However, from the perspective of an observer in $M_L$ near the hypersurface $\mathcal{H}$, locally closed pseudo-timelike loops may appear as the creation of a particle–antiparticle pair at two distinct points $\hat{q},q \in \mathcal{H}$. One can visualize this as an object entering the Riemannian region, temporarily leaving the Lorentzian causal structure, and then re-emerging so as to propagate “forward in time” with respect to the time orientation in $M_L$. Since the mathematical definition of the loop itself does not privilege a direction of traversal, this picture is symmetric. This interpretation highlights that, while pseudo-timelike loops exist mathematically across the signature-changing hypersurface (see Subsection~\ref{Subsec: Local Loops}), they do not correspond to physically realizable violations of causality in $M_L$; instead, they provide a consistent physical picture in which $M_R$ either acts as a region in causal contact with $M_L$ or as a causal barrier.

~
\section*{Acknowledgments}
The authors thank the reviewers for their insightful comments and critique, which have helped to refine and improve this article.
NER is greatly indebted to Richard Schoen for generously welcoming
her into his research group and to Alberto Cattaneo for affording
her with creative independence throughout the duration of this research
endeavor. Moreover, NER acknowledges partial support of the SNF Grant
No. 200021-227719. This research was (partly) supported by the NCCR
SwissMAP, funded by the Swiss National Science Foundation.

~
\textbf{~}\\\textbf{Data availability}: No data was used for the research described in the article.\\
\textbf{~}\\\textbf{Conflict of interest}: The authors have no conflict of interest to declare that are relevant to the content of this article.

~


\end{document}